\documentclass[11pt]{smfart}
\def\pasdegrille{\let\grille = \pasgrille}
\def\ecriture#1#2{\setbox1=\hbox{#1}
\dimen1= \wd1
\dimen2=\ht1
\dimen3=\dp1
\grille #2 \box1 }
\def\aat#1#2#3{
\divide \dimen1 by 48
\dimen3=\dimen1
\multiply \dimen1 by #1
\advance \dimen1 by -\dimen3
\divide \dimen1 by 101
\multiply \dimen1 by 100
\divide \dimen2 by \count11
\multiply \dimen2 by #2 
\setbox0=\hbox{#3}\ht0=0pt\dp0=0pt
  \rlap{\kern\dimen1 \vbox to0pt{\kern-\dimen2\box0\vss}}\dimen1= \wd1
\dimen2=\ht1}
\def\pasgrille{
\count12= \dimen1 
\divide \count12 by 50
\divide \dimen2 by \count12
\count11 =\dimen2
\ 
\divide \dimen1 by 48
\setlength{\unitlength}{\dimen1}
\smash{\rlap{\ }}
\dimen1= \wd1
\dimen2=\ht1
}
\def\grille{
\count12= \dimen1 
\divide \count12 by 50
\divide \dimen2 by \count12
\count11 =\dimen2
\ 
\divide \dimen1 by 48
\setlength{\unitlength}{\dimen1}
\smash{\rlap{\graphpaper[1](0,0)(50, \count11)}}
\dimen1= \wd1
\dimen2=\ht1
}

\pasdegrille
\usepackage{amssymb}
\usepackage{amsmath, amsthm, amsopn, amsfonts}
\usepackage[dvips]{graphicx}
\usepackage{graphpap}
\usepackage[english, francais]{babel}
\newtheorem{theoreme}{Theorem}
\renewcommand{\theequation}{\arabic{equation}}

\newtheorem{proposition}{Proposition}
\newtheorem{lemme}[proposition]{Lemma}

\newtheorem{definition}[proposition]{Definition}
\newtheorem{remarque}[proposition]{Remark}

\numberwithin{equation}{section}
\numberwithin{proposition}{section}
\def\no{\|}

\def\11{{\rm 1~\hspace{-1.4ex}l} }
\def\R{\mathbb R}

\def\Z{\mathbb Z}
\def\N{\mathbb N}

\def\T{\mathbb T}
\def\O{\mathbb O}

\def\cal{\mathcal}

\begin{document}
\selectlanguage{english}
\title[NLS on surfaces]
{ Bilinear eigenfunction estimates and the nonlinear schr\"odinger equation on surfaces 
%version \today
}
\author{N. Burq}
\address{Universit{\'e} Paris Sud,
Math{\'e}matiques,
B{\^a}t 425, 91405 Orsay Cedex}
\email{Nicolas.burq@math.u-psud.fr}
\urladdr{http://www.math.u-psud.fr/~burq}
\author{P. G{\'e}rard}
\address{Universit{\'e} Paris Sud,
Math{\'e}matiques,
B{\^a}t 425, 91405 Orsay Cedex}
\email{Patrick.gerard@math.u-psud.fr}
\author{N. Tzvetkov}
\address{Universit{\'e} Paris Sud,
Math{\'e}matiques,
B{\^a}t 425, 91405 Orsay Cedex}
\email{Nikolay.tzvetkov@math.u-psud.fr}
\urladdr{http://www.math.u-psud.fr/~tzvetkov}
\begin{abstract} 
We study the cubic non linear Schr\"odinger equation (NLS) on compact surfaces.
On the sphere $\mathbb{S}^2$ and more generally on Zoll surfaces, we prove that, 
for $s>1/4$, NLS is uniformly well-posed in $H^s$, which is sharp on the sphere.
The main ingredient in our proof is a sharp bilinear estimate for Laplace spectral projectors on compact surfaces. 
\end{abstract}
\begin{altabstract} 
On \'etudie l'\'equation de Schr\"odinger non lin\'eaire (NLS) sur une surface compacte.
Sur la sph\`ere $\mathbb{S}^2$ et plus g\'en\'eralement sur toute surface de Zoll, 
on d\'emontre que pour $s>1/4$, NLS est uniform\'ement bien pos\'ee dans $H^s$, ce qui est optimal
sur la sph\`ere. Le principal ingr\'edient de notre d\'emonstration est une estimation bilin\'eaire pour 
les projecteurs spectraux du laplacien sur une surface compacte.
\end{altabstract}
\maketitle
\selectlanguage{english}
\tableofcontents

\section{Introduction.}
Let $(M,g)$ be a compact Riemannian manifold of dimension $2$ without boundary.
In this paper we study the Cauchy problem for the cubic nonlinear Schr\"odinger equation posed on $M$,
\begin{equation}\label{nls}
i\partial _tu+\Delta u=|u|^2u
\quad,\quad u(0,x)=u_0(x)\ \ ,
\end{equation}
where the unknown $u=u(t,x)$ is a complex valued function on $\R \times M$,
and $\Delta $ denotes the Laplace operator associated to the metric $g$ on $M$.
It is classical that  smooth solutions of (\ref{nls}) satisfy the following two conservation laws,
\begin{equation}\label{cons}
\int_M |u(t,x)|^2\, dx=\int_M|u_0(x)|^2\, dx\, , \quad
\int _M|\nabla u(t,x)|_g^2\, dx+\frac{1}{2}\int _M|u(t,x)|^4\, dx=E_0\ .
%\left \{ 
%\begin{aligned} 
%&
%\int_M |u(t,x)|^2\, dx=\int_M|u_0(x)|^2\, dx \,
%\\ 
%&\int _M|\nabla u(t,x)|_g^2\, dx+\frac{1}{2}\int _M|u(t,x)|^4\, dx=E_0\ .
%\end{aligned}
%\right .
\end{equation}
As a consequence, the Sobolev $H^1$ norm of the function $u(t,.)$ on $M$ is 
controlled by the $H^1$ norm of the Cauchy data $u_0$. By combining this observation with
energy estimates and a logarithmic Gronwall lemma, Br\'ezis-Gallouet \cite{BrGa80}  proved 
that, for any $u_0\in C^\infty (M)$, the problem (\ref{nls}) has a unique solution $u\in C^\infty (\R \times M)$.
\\

Our purpose is the study of the dynamics of the flow $u_0\mapsto u$ defined by the latter theorem,
in connection with the geometry of the surface $M$. Indeed, because of the infinite speed of propagation for 
equation (\ref{nls}), it is natural to expect the geometry to play an important role even in the small time dynamics. 
Moreover, it is likely that the results will be  sensitive to  the distance on the phase space, which we shall choose to
be the $H^s$ distance, a natural choice in view of (\ref{cons}) and of the properties of the linear Schr\"odinger group.
This leads to the following definition.
\begin{definition} Let $s$ be a real number. We shall say that the Cauchy problem (\ref{nls}) is 
uniformly well-posed in $H^s(M)$ if, for any bounded subset $B$ of $H^s(M)$, there exists $T>0$ such that
the flow map 
\begin{equation}\label{flow}
u_0\in C^\infty (M)\cap B\mapsto u\in C([-T,T],H^s(M))
\end{equation}
is uniformly continuous when the source space is endowed with the $H^s$ norm, and 
when the target space is endowed with
$$\no u\no _{C_TH^s}=\sup _{|t|\leq T}\no u(t)\no _{H^s}.$$
\end{definition}
The notion of uniformly well-posed Cauchy problem has been recently addressed by several authors in the context of 
various nonlinear  evolution equations
(see {\it e.g.} Birnir-Kenig-Ponce-Svanstedt-Vega~\cite{BKPSV}, Gallagher-G\'erard \cite{GaGe01}, Kenig-Ponce-Vega \cite{KePoVe01}, Lebeau \cite{Le01}, 
Burq-G\'erard-Tzvetkov \cite{BuGeTz01-2}, Christ-Colliander-Tao
\cite{ChCoTa01}, Koch-Tzvetkov \cite{KoTz02}). 
Besides the existence of a local continuous flow map on $H^s$, which corresponds to Hadamard's classical 
notion of wellposedness, this notion can be rephrased as  high frequency stability in $H^s$ in the following sense : 
given two bounded sequences of Cauchy data in $H^s$, whose difference converges to $0$ in $H^s$, the difference between 
the corresponding sequences of solutions converges to $0$ in $H^s$ uniformly on a fixed small time interval. In other words, 
if the problem is uniformly well-posed in $H^s$,  the flow acts in small time  on defects of compactness of bounded sequences in $H^s$. 
On the other hand, the lack of uniform well-podsedness corresponds to some nonlinear instability
of the evolution acting on data involving high frequencies.
\\

At this stage, it is worth to observe that, as far as we know,
the only known examples of uniformly well-posed Cauchy problems are  based
on the convergence of the Picard iteration scheme on some convenient Banach space of functions, 
hence are such that the flow map (\ref{flow}) is in
fact Lipschitz continuous. This will be the case for the results in this paper too.
\\

Let us  review the known results about uniform wellposedness for the problem (\ref{nls}) in the literature, starting with the positive 
results. First, in view of the Sobolev embedding theorem in two space dimensions, a simple application of energy methods leads to  
uniform wellposedness in $H^s$ for every $s>1$ on every compact surface $M$. In the case $s=1$, notice that results of 
Vladimirov \cite{Vl84}  and Ogawa-Ozawa \cite{OgOz91} only imply Hadamard wellposedness. 
More recently, in \cite{BuGeTz01-1} (see also~\cite{BuGeTz03-1}), the authors proved  uniform wellposedness in
$H^s$ for every $s>1/2$. The proof is based on the following Strichartz estimates for the linear Schr\"odinger group on $M$,
derived in \cite{BuGeTz01-1} (see also Staffilani-Tataru \cite{StTa02}) : for every function $v_0$ on $M$, for every finite time 
interval $I$,
\begin{equation}\label{strichartz}
\int _I\left (\int _M|{ e}^{it\Delta }v_0(x)|^q\, dx\right )^{p/q}\, dt\quad \leq \quad C(I)\, \no v_0\no _{H^{1/p}(M)}^p\ ,
\end{equation}
for every pair $(p,q)$ satisfying
$$p>2\quad ,\quad \frac {1}{p}+\frac{1}{q}=\frac{1}{2}\quad .$$
Indeed, these estimates allow to gain essentially half a derivative with respect to the Sobolev embedding,
if one agrees to deal with $L^2_tL^\infty _x$ norms rather than with space-time $L^\infty$ norms,
which does not make strong differences in the analysis of (\ref{nls}).
\\
 
On the other hand, using a different method, Bourgain proved in \cite{Bo93}, \cite{Bo93-1} that,
on the torus ${{\mathbb{T}}}^2=\R ^2/{\Z}^2$,  uniform wellposedness holds in
$H^s$ for every $s>0$. This result should be compared to the corresponding one on the Euclidean plane, due to 
Cazenave-Weissler \cite{CaWe90}, and
which is based on Strichartz estimates with no loss of derivatives, namely estimates (\ref{strichartz}) where
the norm in the right hand side is  the $L^2$ norm instead of the $H^{1/p}$ norm. Notice however that Bourgain's approach 
is quite different, since such estimates are not known (or are known to fail) on ${{\mathbb{T}}}^2$.
\\

Let us come to the negative results. The simplest one concerns the torus (in fact it works as well on the one-dimensional torus)
and states that uniform wellposedness fails in $H^s$ for $s<0$ (see \cite{BuGeTz01-2} and Christ-Colliander-Tao \cite{ChCoTa01}). 
More generally, it seems that a suitable adaptation of a recent work of Christ-Colliander-Tao \cite{ChCoTa03} shows that this result 
holds on every compact surface. Therefore the remaining range of regularity to be discussed is the interval $[0,1/2]$ ; moreover, 
on the torus, the picture is complete since uniform wellposedness in $H^s$
is equivalent to positivity of the regularity $s$ --- apart from the critical regularity $s=0$, a very difficult open problem.
\\
\newline
It is interesting to notice that the situation on the standard sphere is quite different. Indeed, in \cite{BuGeTz01-2}, 
strong concentration of some spherical harmonics 
was used to prove that (\ref{nls}) is not uniformly well-posed in $H^s({{\mathbb{S}}}^2)$ for $s<1/4$ 
(see also the recent work of Banica \cite{Ba03}, which gives more
precise results). The main purpose of this paper is to complete the picture on the sphere.
\begin{theoreme}\label{th1}
If $(M,g)$ is the standard sphere, or more generally  a Zoll surface, then the Cauchy problem for (\ref{nls}) 
is  uniformly well-posed in $H^s(M)$ for every $s>1/4$.
\end{theoreme}
We recall that a Zoll surface is a surface on which the geodesic flow is periodic (see Besse \cite{Be78} for a detailed exposition).
It is worth noticing that the construction on ${{\mathbb{S}}}^2$ of \cite{BuGeTz01-2} can be extended to many others revolution
surfaces and thus Theorem \ref{th1} provides sharp uniform wellposedness results far a large class of Zoll surfaces.
The details of this construction will be given elsewhere (see \cite{BGT-new}). 
\begin{remarque}
In the case of the torus and of the sphere, the above results display a critical regularity threshold
$s_c$ above which  uniform wellposedness holds, and below which uniform wellposedness fails : we have $s_c({{\mathbb{T}}}^2)=0$ and,
in view of Theorem 1, $s_c({{\mathbb{S}}}^2)=1/4$. We expect that such a threshold exists on every surface. Of course the most interesting 
open question is to relate the value of this threshold to the geometric properties of the surface. In particular, it would 
be interesting to know whether there exist surfaces for which this threshold exceeds $1/4$.
\end{remarque}
Let us describe briefly the main steps in the proof of Theorem 1. A first step is to reduce, on every compact surface, 
the uniform wellposedness for (\ref{nls}) in  $H^s$ for every $s>s_0$, to the following bilinear inequality on solutions 
of the linear Schr\"odinger equation,
\begin{equation}\label{bil}
\left (\int _{[0,1]\times M}|{ e}^{it\Delta }f\, (x)\, { e}^{it\Delta }g\, (x)|^2\, dt\, dx\right )^{1/2}
\leq  C\, (\min (N,L))^{s_0}\no f\no _{L^2}\, \no g\no _{L^2}\, ,
\end{equation}
where $N,L$ are large dyadic numbers, and $f,g$ are supposed to be spectrally localized on dyadic intervals of order $N,L$ respectively, 
namely
\begin{eqnarray}\label{lokalizatzia}
\11_{N\leq \sqrt {-\Delta }\leq 2N} (f)=f\quad ,\quad \11_{L\leq \sqrt {-\Delta }\leq 2L} (g)=g\quad .
\end{eqnarray}
Notice that such kind of bilinear estimates were established and used by several authors in the context of the wave equation and of the
Schr\"odinger equation with constant coefficients (see Klainerman-Machedon-Bourgain-Tataru \cite{KlMa96}, Bourgain \cite{Bo93, Bo93-1}, Foschi-Klainerman~\cite{FoKl}, 
Tao \cite{Ta01} and references therein). 
The general principle of the above reduction is based on the systematic use of conormal spaces introduced by 
Bourgain in \cite{Bo93} (see also Ginibre \cite{Gi95}). 
\\

The second step in the proof of Theorem 1 is to observe that, due to the good localization of the spectrum of a Zoll surface 
(see Guillemin \cite{Gu77}, Colin de Verdi\`ere \cite{CdV79}), the above bilinear inequality for $s_0>1/4$ reduces on a Zoll 
surface to the following result on the spectral projectors.
\begin{theoreme}\label{th2}  
Let $(M,g)$ be a compact surface, and let $\chi \in C^\infty _0(\R )$. For each $\lambda  \geq 0$, we introduce the operator
$$\chi _\lambda =\chi (\sqrt {-\Delta }-\lambda)\ .$$
There exist $C>0$ such that, for all $\lambda ,\mu \geq 1$, for all functions $f,g$ on $M$,
\begin{equation}\label{bilsogge}
\no \chi _\lambda f\, \chi _\mu g\no _{L^2(M)}\, \leq
\, C\, (\min (\lambda ,\mu ))^{1/4}\, \no f\no _{L^2(M)}\, \no g\no _{L^2(M)}\quad .
\end{equation}
\end{theoreme}
Notice that the case $\lambda =\mu $ in Theorem \ref{th2} is a particular case of general 
$L^p$ estimates due to Sogge (\cite{So86, So88, So93})
which take the following form, for every $p\in [2,+\infty ]$,
\begin{equation}\label{sogge}
\no \chi _\lambda f\no _{L^2(M)}\, \leq\, C\, \lambda ^{s(p)}\, \no f\no _{L^2(M)}\, ,\, \lambda \geq 1,
\end{equation}
where $s(4)=1/8$, and more generally $s(p)$ is given in terms of $1/p$ by the following diagram:
\begin{figure}[ht]
\label{fig:burq}
$$\ecriture{\includegraphics[width=7cm]{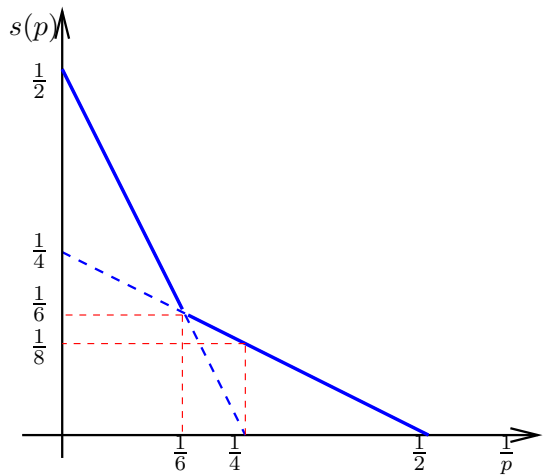}}
{\aat{2}{35}{$\frac 1 2$}\aat{2}{19}{$\frac 1 4$}\aat{2}{14}{$\frac 1 6$}\aat{2}{10}{$\frac 1 8$}\aat{15}{0}{$\frac 1 6$}\aat{20}{0}{$\frac 1 4$}\aat{37}{0}{$\frac 1 2$}\aat{0}{40}{$s(p)$}\aat{45}{0}{$\frac 1 p$}}$$
\caption{The Sogge diagram in dimension $2$}
\end{figure}
 
Of course the advantage of using estimate (\ref{bilsogge}), rather than applying H\"older's inequality and (\ref{sogge}), is to
make the highest frequency disappear from the right hand side, which is crucial in the nonlinear analysis of the first step.
A glance at the particular case of the sphere is particularly enlightening. Indeed, in this case, estimate (\ref{bilsogge})
means that the $L^2$ norm of the product of two spherical harmonics divided by the product of their $L^2$ 
norms is bounded by $l^{1/4}$, where $l$ is the minimum of their two degrees. It is interesting to notice that
estimate (\ref{bilsogge}) (as well as (\ref{sogge})) is made optimal by choosing the following spherical harmonics,
\begin{equation}\label{quasimode}
\phi _n(x)=(x_1+ix_2)^n \, ,\, x=(x_1,x_2,x_3)\in \R ^3 ,\, x_1^2+x_2^2+x_3^2=1\, ,
\end{equation}
which concentrates onto the large circle $\{ x_3=0\}$ of ${{\mathbb{S}}}^2$.
Indeed, this kind of spherical harmonics already appeared in a paper by Stanton and Weinstein \cite{StWe81} which was
a motivation for Sogge's result \cite{So86}. Inspired by these earlier works, we established the nonlinear instability
in \cite{BuGeTz01-2} precisely by studying the Cauchy problem (\ref{nls}) with $u_0=\phi _n$. Therefore, in some sense, the behavior of
$\phi _n$ can summarize in itself the whole problem of uniform wellposedness on the sphere. It is also interesting to notice that
the same example shows that a bilinear inequality such as (\ref{bilsogge}) with an $L^p$ norm in the left hand side
could not hold if $p>2$. In this case, the high frequency would necessarily appear in the right hand side.
\\

The third and main step in the proof of Theorem 1 is of course  Theorem 2 itself. By a reduction which is now classical
(see \cite{So93}), one is led to studying oscillatory integral operators of Carleson-Sj\"olin's type (see Carleson-Sj\"olin \cite{CaSj72},
H\"ormander \cite{Ho73}). At this stage our proof consists in applying the 
calculus of Fourier integral operators (see H\"ormander \cite{Ho71}) in order to reduce the problem to a bilinear 
Carleson-Sj\"olin inequality, which we then prove by adapting an idea used by H\"ormander in \cite{Ho73}.
\\

This paper is organized as follows. In section 2, we precise the relationship between uniform wellposedness
and bilinear inequalities (\ref{bil}) on solutions of the linear Schr\"odinger equation. In section 3, we take advantage
of the localization of the spectrum of a Zoll surface in order to deduce Theorem 1 from Theorem 2. Finally, section 4 is devoted to the proof of Theorem 2.
We close this introduction by mentioning that the results of this paper were announced in \cite{BuGeTz03}.
\\

{\bf Acknowledgements.} We are grateful to Terence Tao for valuable conversations about bilinear estimates. 

\section{The non linear analysis}
In this section $M$ is any compact surface . 
However the analysis can be easily adapted to more general situations as higher dimensional 
manifolds, focusing nonlinearities, etc.
We prove that the uniform wellposedness of ~\eqref{nls} on $M$ can be 
deduced from bilinear estimates on the solutions of the linear equation. This is a result of independent interest that we state below: 
\begin{theoreme}\label{th3} 
Suppose that there exists $C>0$ and $s_0$ such that for any $u_0, v_0 \in L^2( M)$ satisfying
$$
\11_{\sqrt{-\Delta}\in [N,2N]}(u_0)=u_0,\quad
\11_{\sqrt{-\Delta}\in [L,2L]}(v_0)=v_0
$$ 
one has (with $s_{0}\geq 0$)
\begin{equation}\label{eq.bilin}
\|e^{it \Delta}u_0\,\, e^{it\Delta} v_0\|_{L^2((0,1)_t\times M)}
\leq 
C(\min(N,L))^{s_0}\|u_{0}\|_{L^2(M)} \|v_{0}\|_{L^2(M)}.
\end{equation}
Then the Cauchy problem~\eqref{nls} is uniformly well-posed in $H^s( M)$ for any $s>s_0$.
\end{theoreme}
The proof of Theorem~\ref{th3} relies on the use of Bourgain spaces. We first show that
~\eqref{eq.bilin} is equivalent to a bilinear estimate in Bourgain spaces.
We then prove the crucial non linear estimate which yields the proof of the uniform wellposedness
due to a contraction argument in a suitable Bourgain space associated to $H^s$ regularity.
\subsection{Bourgain spaces}
Let $A$ be a non negative self adjoint operator with compact resolvent on $L^2(M)$. 
Denote by $(e_{k})$ an $L^2$ orthonormal basis of eigenfunctions of $A$ associated to eigenvalues $\mu_k$, 
by $P_k$ the orthogonal projector on $e_{k}$, and by $H_{A}^s(M)$ the natural Sobolev space associated to $(\text{Id}+A)^{1/2}$,
equipped with the norm
$$\|u\|_{H^s_A(M)}^2= \sum_k \langle \mu_{k}\rangle^s \|P_k u\|_{L^2(M)}^2.$$ 
\begin{definition}
 The space $X^{s,b}_{A}(\R\times M)$ is the completion of 
$C^\infty_{0}({\mathbb R}_{t}; H^s_{A}(M))$  for the  norm
\begin{equation}
\label{eq2.200}
\begin{aligned}
\|u\|_{{X}^{s,b}_{A}(\R\times M)}^2 &
=
\sum_{k}\|\langle\tau+\mu_{k}\rangle^{b}\langle \mu_{k}\rangle^{s }\widehat{P_k u}(\tau)\|^2_{L^{2}(\R_{\tau}\,;\,L^{2}(M))}
\\
& = \|e^{it A} u(t, \cdot)\|^2_{H^b( {\mathbb R}_{t}\,;\,H^s_{A}(M))},
\end{aligned}
\end{equation}
where $\widehat{P_k u}(\tau)$ denotes the Fourier transform of $P_k u$ with respect to the time variable.
\end{definition} 
\begin{remarque}\label{re2.1} The definition (and the norm) of the space $X^{s,b}_{A}$ clearly depends on the operator $A$. 
However if $A,B$ are two operators as above, having the same eigenfunctions and such that the eigenvalues $\lambda_k$ and $\mu_k$  
of $A$ and $B$ respectively satisfy 
$$\exists\, C>0\,:\, \forall k,\, |\lambda_k- \mu_k| \leq C$$
then
$$
\exists\, C>0\,:\, 
\forall k \in {\mathbb N}, \,\forall \tau \in {\mathbb R},\, 
\frac 1 C \langle\tau + \lambda_{k}\rangle \leq \langle \tau + \mu_{k}\rangle \leq C \langle \tau + \lambda_{k} \rangle 
$$
and consequently $X^{s,b}_{A}=X^{s,b}_{B}$ (with equivalent norms). 
\end{remarque}
We next reformulate the bilinear estimate~\eqref{eq.bilin} in the $X^{s,b}$ context.
\begin{lemme}\label{le2.5}
Let $s\in\R$.
The two following statements are equivalent :
\begin{itemize} 
\item 
For any $u_0, v_0 \in L^2( M)$ satisfying
$$
\11_{\sqrt{A}\in [N,2N]}(u_0)=u_0,\quad
\11_{\sqrt{A}\in [L,2L]}(v_0)=v_0
$$ 
one has 
\begin{equation}
\label{eq2.1}
\|e^{-itA} u_0\,\, e^{-itA} v_0\|_{L^2((0,1)_t\times M)}\leq C(\min(N,L))^{s}\|u_0\|_{L^2}\|v_0\|_{L^2}.
\end{equation}
\item 
For any $b>1/2$ and any $f, g \in X^{0,b}(\R\times M)$ satisfying
$$
\11_{\sqrt{A}\in [N,2N]}(f)=f,\quad \11_{\sqrt{A}\in [L,2L]}(g)=g
$$ 
one has 
\begin{equation}\label{eq2.2}
\|f\, g\|_{L^2(\mathbb{R}\times M)}\leq C(\min(N,L))^{s}
\|f\|_{X^{0,b}_{A}(\R\times M)}\|g\|_{X^{0,b}_{A}(\R\times M)}.
\end{equation}
\end{itemize}
\end{lemme}
\proof 
If $u(t)=e^{itA} u_0$ then for any $\psi \in C^\infty_0( \mathbb{R})$ and any $b$, $\psi(t) u(t)\in X^{0,b}_{A}(\R_{t}\times M)$ 
with 
$$\|\psi\, u\|_{X^{0,b}_{A}(\R\times M)}\leq C \|u_0\|_{L^2(M)},$$
which shows that~\eqref{eq2.2} implies \eqref{eq2.1}. 
To show the reverse implication, suppose first that $f(t)$ and $g(t)$ are supported in time in the interval $(0,1)$ and write
$$
f(t)= e^{-itA}e^{itA} f(t):= e^{-itA}F(t),\quad g(t)= e^{-itA}e^{itA} g(t):= e^{-itA}G(t)
$$
Then 
$$
f(t)= \frac{1}{2\pi}\int_{-\infty}^{\infty} e^{it \tau} e^{-itA} \widehat{F} ( \tau) d\tau,\quad 
g(t)= \frac{1}{2\pi}\int_{-\infty}^{\infty} e^{it \tau} e^{-itA} \widehat{G} ( \tau) d\tau
$$
and hence
$$(f\,g)(t) = \frac{1}{(2\pi)^2}
\int_{-\infty}^{\infty}\int_{-\infty}^{\infty}
e^{it(\tau + \sigma)}\, e^{-it A }\widehat{ F} (\tau)\,\, e^{-itA} \widehat{G}( \sigma) d\tau d\sigma. 
$$
Ignoring the oscillating factor $e^{it(\tau + \sigma)}$,
using \eqref{eq2.1} and the Cauchy-Schwarz inequality in $(\tau,\sigma)$ (in this place we use that $b>1/2$ to get the
needed integrability) yields
\begin{equation}
\begin{aligned} 
\|fg\|_{L^2((0,1)_t\times M)} 
&\leq 
C(\min(N,L))^s \int_{\tau , \sigma}\|\widehat{F} (\tau)\|_{L^2(M)} \|\widehat{G}( \sigma)\|_{L^2(M)} d\tau d\sigma\\
&\leq C (\min(N,L))^{s} \|\langle\tau\rangle^{b}\widehat{F} (\tau)\|_{L^2( {\mathbb R}_\tau \times M)}
\|\langle\sigma\rangle^b\widehat{G}( \sigma)\|_{L^2(\mathbb{R}_\sigma \times M)}\\
&= C(\min(N,L))^{s} \|f\|_{X^{0,b}_{A}(\R\times M)} \|g\|_{X^{0,b}_{A}(\R\times M)}.
\end{aligned}\end{equation}
Finally, by decomposing $f(t)= \sum_{n\in \mathbb{Z}} \psi (t-n/2)f(t)$ and 
$g(t)=\sum_{n\in \mathbb{Z}}\psi(t-n/2)g(t)$ with a suitable $\psi\in C^\infty_0( \mathbb{R})$ 
supported in $(0,1)$, the general case for $f(t)$ and $g(t)$ follows from the considered particular case
of $f(t)$ and $g(t)$ supported in time in the interval $(0,1)$.
\qed
\subsection{The crucial nonlinear estimate in $X^{s,b}$}
In this section we come back to the proof of Theorem \ref{th3}.
Therefore we consider the Bourgain spaces associated to $A=-\Delta$ and we drop the subscript $A$.
We begin with an elementary observation.
\begin{lemme}\label{prop2.1}
For every $b>\frac{1}{4}$, there exists $C>0$ such that for every
$u\in X^{0,b}(\R\times M)$,
\begin{equation}\label{eq2.100}
\|u\|_{L^{4}(\R\,;\,L^{2}(M))}
\leq C
\|u\|_{X^{0,b}(\R\times M)}.
\end{equation}
\end{lemme}
\proof
Write the Fourier transform inversion as follows
$$
P_k u(t)=
\frac{1}{2\pi}
\int_{-\infty}^{\infty}
\frac{\langle\tau+\mu_k\rangle^{b}}{\langle\tau+\mu_k\rangle^{b}}
\widehat{P_k u}(\tau)e^{it\tau} d\tau.
$$
Therefore for $b>\frac{1}{2}$, we get by the Cauchy-Schwarz inequality, applied in $\tau$,
\begin{equation}\label{eq2.400}
|P_k u(t)|\lesssim\left\{\int_{-\infty}^{\infty}
\langle\tau+\mu_k\rangle^{2b}|\widehat{P_k u}(\tau)|^{2}d\tau\right\}^{\frac{1}{2}}.
\end{equation}
Squaring (\ref{eq2.400}), integrating over $M$ and summing over $k$ yields,
\begin{equation}\label{eq2.3}
\|u\|_{L^{\infty}(\R\,;\,L^{2}(M))}\lesssim\|u\|_{X^{0,b}(\R\times M)},\quad b>\frac{1}{2}.
\end{equation}
We also have the trivial identity,
\begin{equation}\label{eq2.4}
\|u\|_{L^{2}(\R\,;\,L^{2}(M))}=\|u\|_{X^{0,0}(\R\times M)}.
\end{equation}
Interpolation between (\ref{eq2.3}) and (\ref{eq2.4}) gives (\ref{eq2.100}).
\qed
\\

Let us now state the main result of this section.
\begin{proposition}\label{prop2.2}
Under the assumptions of Theorem \ref{th3},
let $s>s_{0}$. There exists $(b,b')\in \R^2$ satisfying
\begin{equation}\label{star}
0<b'<\frac{1}{2}<b,\quad b+b'<1,
\end{equation}
and $C>0$ such that for every triple $(u_j)$, $j=1,2,3$, in $X^{s,b}(\R\times M)$,
\begin{equation}\label{eq2.5}
\|u_1\, \overline{u_2}\, u_3\|_{X^{s,-b'}} \leq C\|u_1\|_{X^{s,b}}\|u_2\|_{X^{s,b}}\|u_3\|_{X^{s,b}}.
\end{equation}
\end{proposition}
\proof Define $w_j$ as follows
$$
w_{j}(t)=
\frac{1}{2\pi}
\sum_{k\in\mathbb{N}}\,
\langle \mu_{k}\rangle^{s/2}
\int_{-\infty}^{\infty}
\langle\tau+\mu_{k}\rangle^{b}
\widehat{P_k u_j}(\tau)e^{it\tau}d\tau,\quad j=1,2,3
$$
or equivalently
$$
u_{j}(t)=
\frac{1}{2\pi}
\sum_{k\in\mathbb{N}}\,
\langle \mu_{k}\rangle^{-s/2}
\int_{-\infty}^{\infty}
\langle\tau+\mu_{k}\rangle^{-b}
\widehat{P_k w_j}(\tau)e^{it\tau}d\tau,\quad j=1,2,3.
$$
Notice that
$$
\|u_j\|_{X^{s,b}}=
\|w_j\|_{L^2},\quad j=1,2,3.
$$
A duality argument reduces (\ref{eq2.5}) to
\begin{equation}\label{eq2.6}
\left|\int_{\R\times M}\overline{u_0}\,u_1\,\overline{u_2}\,u_3\right|\lesssim\|w_0\|_{L^{2}(\R\times M)}
\prod_{j=1}^{3}
\|w_j\|_{L^{2}(\R\times M)},
\end{equation}
where $w_0$ is arbitrary in $L^{2}(\R\times M)$ and
$$
u_{0}(t)=
\frac{1}{2\pi}
\sum_{{k}\in\mathbb{N}}\,
\langle \mu_{k}\rangle^{s/2}
\int_{-\infty}^{\infty}
\langle\tau+\mu_{k}\rangle^{-b'}
\widehat{P_k w_0}(\tau)e^{it\tau}d\tau.
$$
The next step is to perform a localization in the integral of the left hand-side of  (\ref{eq2.6}). 
In the sequel $N_0$, $N_1$, $N_2$, $N_3$, $L_0$,  $L_1$, $L_2$, $L_3$ will denote dyadic integers, i.e.
$N_j=2^{n_j}$, $n_j\in\N$, $L_j=2^{l_j}$, $l_j\in\N$, $j=0,1,2,3$.
We set
$$
N:=(N_0,N_1,N_2,N_3),\quad
L:=(L_0,L_1,L_2,L_3)
$$
and that notation will be frequently used in the rest of the proof.
By summation over $N$, we mean summation over all possible dyadic values of $N_0$, $N_1$, $N_2$, $N_3$.
Similar convention will be adopted for the summation over $L$. Denote by $J$ the left hand-side of  (\ref{eq2.6}).
Then we can write
$$
J\lesssim
\sum_{L}\sum_{N}I(L,N),
$$
where
\begin{eqnarray}\label{iii}
I(L,N)=
\left|\int_{\R\times M}
\overline{u_0^{L_0 N_0}}\,
u_1^{L_1 N_1}\,
\overline{u_2^{L_2 N_2}}\,
u_3^{L_3 N_3}
\right|\
\end{eqnarray}
with
\begin{eqnarray}\label{a}
\qquad
u_{j}^{L_j N_j}(t)=
\frac{1}{2\pi}
\sum_{N_j\leq \langle\mu_k\rangle^{1/2}< 2N_j}\,
\langle \mu_k\rangle^{-s}
\int_{L_j\leq \langle\tau+\mu_k\rangle< 2L_j}
\langle\tau+\mu_k\rangle^{-b}
\widehat{P_k w_j}(\tau)e^{it\tau}d\tau
\end{eqnarray}
for $j=1,2,3$ and
\begin{eqnarray}\label{b}
\qquad
u_{0}^{L_0 N_0}(t)=
\frac{1}{2\pi}
\sum_{N_0\leq \langle\mu_k\rangle^{1/2}< 2N_0}\,
\langle \mu_k\rangle^{s}
\int_{L_0\leq \langle\tau+\mu_k\rangle< 2L_0}
\langle\tau+\mu_k\rangle^{-b'}
\widehat{P_k w_0}(\tau)e^{it\tau}d\tau.
\end{eqnarray}
Notice that
$$
\|u_{j}^{L_j N_j}\|_{X^{\sigma,\beta}}
\lesssim
L_{j}^{\beta-b}N_{j}^{\sigma-s}c_{j}(L_j,N_j),\quad j=1,2,3
$$
and
$$
\|u_{0}^{L_0 N_0}\|_{X^{\sigma,\beta}}
\lesssim
L_{0}^{\beta-b'}N_{0}^{\sigma+s}c_{0}(L_0,N_0),
$$
where the sequences of real numbers $c_{j}(L_j,N_j)$ satisfy
$$
\sum_{L_j}\sum_{N_j}
(c_{j}(L_j,N_j))^{2}
\lesssim
\|w_j\|_{L^2}^{2},\quad j=0,1,2,3.
$$
In particular
$c_{j}(L_j,N_j)\lesssim\|w_j\|_{L^2},$ $j=0,1,2,3.$
\\

Our next step is to give a bound for the product of four functions on $M$, 
localized in frequency on four intervals $I_k$ and such that one of
the intervals is much shifted with respect to the others.
\begin{lemme}\label{le2.1}
There exists $C>0$ such that, if for any $j=1, 2, 3$, $ C\mu_{k_{j}}\leq \mu_{k_{0}}$, then for every $p>0$ 
there exists $C_{p}>0$ such that for every
$w_j\in L^{2}(M)$, $j=0,1,2,3$,
\begin{equation}\label{eq2.10}
\left|\int_{M}P_{k_0}w_0\,P_{k_1}w_1\,P_{k_2}w_2\,P_{k_3}w_3\, dx \right|
\leq
C_{p}\,
\mu_{k_{0}}^{-p}\prod_{j=0}^{3}\|w_j\|_{L^2}.
\end{equation}
\end{lemme}
\proof
There are several ways to prove Lemma \ref{le2.1}. Our proof relies on the
following result of  Sogge~\cite[Chap. 5.1]{So93}, which we shall use in its full
strength in Section \ref{sec4}.
\begin{proposition}\label{propsogge}
There exists a function $\chi \in \mathcal{S}( \mathbb{R})$ equal to $1$ at the point $0$ such that 
$$
\chi (\sqrt{-\Delta}-\lambda)f=\lambda^{\frac{1}{2}}{\mathcal {T}}_{\lambda}f+R_{\lambda}f,
$$
with
$$
\forall\, p,\, s \in \mathbb{N}\,,\,\, \exists\, C_{p,s} >0\,:\,  
\|R_{\lambda}f\|_{H^s(M)}\leq C_{p,s} \lambda^{-p}\|f\|_{L^2}
$$
and in a coordinate system close to $x_0\in M$,
$$
{\mathcal {T}}_{\lambda}f(x)=\int_{\R^2}e^{i\lambda\varphi(x,y)}a(x,y,\lambda)f(y)dy,
$$
with the following asymptotic expansion in $C^\infty_0$, 
$$a(x,y,\lambda)\sim \sum_{j\geq 0} a_{j}(x,y)\lambda ^{-j},$$
with
$
a_j\in C_{0}^{\infty}
$ 
supported in
$$
\frac{\varepsilon}{C}
\leq |x-y|
\leq 
C\varepsilon,
$$
and $-\varphi(x,y)=d_{g}(x,y)$ is the geodesic distance between $x$ and $y$.
\end{proposition}
We fix $0<\varepsilon<1/2$ and distinguish now different cases in (\ref{eq2.10})  according whether 
$ \mu_{k_{0}}^\varepsilon< C\mu_{k_{j}}<\mu_{ k_{0}}$ or $\mu_{k_{j}}< \mu_{k_{0}}^ \varepsilon$. 
In the first case, using that 
$$
\chi(\sqrt{ - \Delta}- \mu_{k_{j}}^{1/2})P_{k_{j}}= P_{k_{j}},
$$  
we can replace $P_{k_{j}}w_{j}$ by ${\mathcal {T}}_{\mu_{k_{j}}^{1/2}}\widetilde{w}_{j}$ with  $\widetilde{w}_{j}=P_{k_{j}}w_{j}$, 
modulo the error term $R_{\mu_{k_{j}}^{1/2}}$ which gives for any $p$ an error of order $\mathcal{O}(\mu_{k_{j}}^{-p})$ and 
consequently of order $\mathcal{O}({\mu_{k_{0}}^{-p}})$ for any $p$. 
In the second case we keep $P_{k_{j}}w_{j}$. To fix ideas suppose for example that 
$\mu_{k_{0}}^ \varepsilon < C\mu_{k_{1}}<  \mu_{k_{0}}$ and $\mu_{k_{2}}, \mu_{k_{3}}<\mu_{k_{0}}^{\varepsilon}$ 
(the other cases being similar). Using a partition of unity, we have consequently to estimate:
\begin{multline}\label{eq2.11}
\left|\int_{x\in U}\prod_{j=0}^{1}(\mathcal{T}_{\mu_{k_{j}}^{1/2}}\widetilde{w}_j)(x) \prod_{j=2}^3 P_{k_{j}}w_{j}(x) dx\right|=
\\
\left|
\int_{y_0}\int_{y_1}
\int_{x\in U}
e^{i\mu_{k_{0}}^{1/2}\Phi(x,y_{0}, y_{1})}\left(
\prod_{j=0}^{1}a(x,y_j,\mu_{k_j}^{1/2})\widetilde{w}_{j}(y_j)dy_j 
\right) P_{k_{2}}w_{2}(x)P_{k_{3}}w_{3}(x) 
dx
\right|,
\end{multline}
with
$$
\Phi=
\varphi(x,y_0)+\sqrt{\frac{\mu_{k_{1}}}{ \mu_{k_{0}}} }\varphi(x,y_1).
$$
and $a_{k_{j}}$ supported in $U$ (in the $x$ variable).
\\
\newline
Since $|\nabla_{x}\varphi|$ is uniformly bounded from below and bounded from above together with all its derivatives, we obtain that
for $C>0$ large enough, there exist $c_1>0$ and $C_{\beta}$ such that
$$
|\nabla_x \Phi|\geq c_1,\quad |\partial^{\beta}_{x}\Phi|\leq C_{\beta}.
$$
It remains to perform integrations by parts in the $x$ variable in~\eqref{eq2.11}. 
Each such integration  gains a power of $\mu_{k_{0}}^{-1/2}$ and looses 
(due to the derivatives on $\prod_{j=2}^3 P_{k_{j}}w_{j}(x)$) a power of 
$\sqrt{\text{max}( \mu_{k_{2}}, \mu_{k_{3}})}\leq \mu_{k_{0}}^{\varepsilon/2}$. 
After a suitable number of integration by parts, we can conclude using a crude Sobolev embedding.\qed
\\

Lemma \ref{le2.1} is now the key for the proof of the next statement.
\begin{lemme}\label{le2.2}
There exists $C>0$ such that if $N_0\geq C(N_1+N_2+N_3)$ then for every $p>0$ there exists $C_{p}$ such that
$$
I(L,N)\leq C_{p}N_{0}^{-p}\frac{(L_1 L_2)^{\frac{1}{2}}}{L_0^{b'} (L_1 L_2 L_3)^{b}}
\prod_{j=0}^{3}\|w_j\|_{L^2}.
$$
\end{lemme}
\begin{remarque}
We note that if $M=\mathbb{S}^2$ then Lemma \ref{le2.2} is trivial since in that case, 
by an elementary observation on the degree of the corresponding spherical harmonics, 
we obtain that if $N_0> 2(N_1+N_2+N_3)$ then $I(L,N)=0$.
\end{remarque}
\proof
We substitute~\eqref{a}, \eqref{b} into ~\eqref{iii}. We obtain
$$I(L, N) \lesssim \left(\frac {N_0} {N_1 N_2 N_3}\right) ^s  \sum _{(k_0, k_1, k_2, k_3) \in S(L, N)} \widetilde {I} ( k_0, k_1, k_2, k_3),
$$
where 
$$ S(L, N)= \{ (k_0, k_1, k_2, k_3) \in \mathbb{N} ^4 \, :\, N_j \leq \langle \mu_{k_j} \rangle ^{1/2} \leq 2 N_j; \, j=0, \dots, 3\}
$$
and 
\begin{equation} \widetilde{I} (k_0, k_1, k_2, k_3) = \left| \int_D{
\int_M P_{k_0} \overline{ \widehat{ w}_{k_0}(\tau_0)} P_{k_1}{ \widehat{ w}_{k_1}(\tau_1)} P_{k_2}\overline{ \widehat{ w}_{k_2}(\tau_2)} P_{k_3}{ \widehat{ w}_{k_3}(\tau_3)}}dx d\mu \right|
\end{equation}
where the integral in $(\tau_0,\tau_1, \tau_2,\tau_3) $ is over the domain
$$
D= \{
(\tau_0,\tau_1, \tau_2,\tau_3)\,:\, L_j \leq \langle \tau_j + \mu_{k_j}\rangle \leq 2 L_j, j=0,\cdots, 3
\}
$$
endowed with the measure
$$ d\mu =
 \frac{ \delta ( - \tau_0+ \tau_1 -\tau_2+\tau_3) d\tau_0 d\tau_1 d\tau_2d\tau_3} {{ \langle \tau_0+ \mu_{k_0}\rangle^{b'} \prod_{i=1}^3 \langle \tau_i+ \mu_{k_i}\rangle^{b}}}.
$$
Next we apply Lemma \ref{le2.1} and we introduce functions $h_{j}(\tau)=\|\widehat{w_j}(\tau)\|_{L^{2}(M)}$, $j=0,1,2,3$. We obtain
\begin{equation*}
\widetilde{I} (k_0, k_1, k_2, k_3)\leq 
\frac{C_{p}N_{0}^{-p}}{L_0^{b'} (L_1 L_2 L_3)^{b}}
\int_{\Lambda}
h_{0}(\tau)h_{1}(\tau-\tau_1)h_{2}(\tau_2-\tau_1)h_{3}(\tau_2)d\tau d\tau_1 d\tau_2\,,
\end{equation*}
where 
$$
\Lambda=\{
(\tau,\tau_1,\tau_2)\,:\,
L_1\leq \langle\tau-\tau_1+\mu_{k_1}\rangle\leq 2L_1,\,\,
L_2\leq \langle\tau_2-\tau_1+\mu_{k_2}\rangle\leq 2L_2
\}
$$
(note that we simply neglected the localizations with respect to $L_0$ and $L_3$).
\\
Next we evaluate the integral over $\tau_1$, $\tau_2$, $\tau$.
We first apply the Cauchy-Schwarz inequality with respect to $(\tau_1,\tau_2)$ and then another Cauchy-Schwarz inequality
with respect to $\tau$ to obtain
\begin{eqnarray*}
\int_{\Lambda}
h_{0}(\tau)h_{1}(\tau-\tau_1)h_{2}(\tau_2-\tau_1)h_{3}(\tau_2)d\tau d\tau_1 d\tau_2
\leq
(\sup_{\tau}\alpha(\tau))^{1/2}
\prod_{j=0}^{3}\|h_{j}\|_{L^{2}(\R)},
\end{eqnarray*}
where
$$
\alpha(\tau)=
{\rm measure}\{(\tau_1,\tau_2)\in \R^2\,:\, L_1\leq \langle\tau-\tau_1+\mu_{k_1}\rangle\leq 2L_1,\,\,
L_2\leq \langle\tau_2-\tau_1+\mu_{k_2}\rangle\leq 2L_2\}.
$$
Clearly $\alpha(\tau)\leq L_1 L_2$ and therefore
\begin{equation*}
\widetilde{I} (k_0, k_1, k_2, k_3)
\leq
\frac{C_{p}N_{0}^{-p}(L_1 L_2)^{1/2}}{L_0^{b'} (L_1 L_2 L_3)^{b}}
\prod_{j=0}^{3}\|w_{j}\|_{L^{2}}.
\end{equation*} Since the Weyl asymptotics gives 
$$
\# \{k\,:\,N_j\leq \langle\mu_k\rangle^{1/2}\leq 2N_j\}\lesssim N_j^{2}\, ,
$$
the number of elements in the set $S(L,N)$ is bounded by $(N_0 N_1 N_2 N_3)^2$, which completes the proof of Lemma~\ref{le2.2}.
\qed
\\
 
In view of lemma \ref{le2.2}, it is natural to write
$$
J\leq J_1 +J_2,
$$
where the $N$ summation is restricted to $N_0\geq C(N_1+N_2+N_3)$ in $J_1$ and all other possibilities are in $J_2$.
If $b>\frac{1}{2}$ and $b'>0$, then lemma \ref{le2.2} enables one to bound $J_1$ by
$
\prod_{j=0}^{3}\|w_j\|_{L^2}
$
using a simple summation of geometric series in all dyadic variables.
Hence it remains to bound
$$
J_{2}=
\sum_{L}\quad
\sum_{N_0\leq C(N_1+N_2+N_3)}
I(L,N),
$$
The conjugates being of no importance in the proof below, we can use a symmetry argument and suppose that the summation in $J_{2}$ is restricted to $N_1\geq N_2\geq N_3$.
Our next step is to perform two different ways of evaluating $I(L,N)$. 
The first one is better with respect to the $N$ localization.
The second one gives a better result on the $L$ localization.
The useful one will be a suitable interpolation between the two bounds.
\\

For $N_1\geq N_2\geq N_3$, using the Cauchy-Schwarz inequality, the assumption in Theorem~\ref{th3}, 
and Lemma~\ref{le2.5}, we obtain, for every $\epsilon>0$,
\begin{eqnarray*}
I(L,N) 
& \leq & 
\|u_{0}^{L_0 N_0}u_{2}^{L_2 N_2}\|_{L^{2}(\R\times M)}\|u_{1}^{L_1 N_1}u_{3}^{L_3 N_3}\|_{L^{2}(\R\times M)}
\\
& \lesssim &
C_{\epsilon}
\frac{(L_0 L_1 L_2 L_3)^{\frac{1}{2}+\epsilon}}
{L_0^{b'}(L_1 L_2 L_3)^{b}}
\frac{N_0^s}{(N_1 N_2 N_3)^{s}}
(N_2 N_3)^{s_0}
\prod_{j=0}^{3}c_{j}(L_j,N_j).
\end{eqnarray*}
On the other hand, for $N_1\geq N_2\geq N_3$, using the H\"older inequality, lemma \ref{prop2.1} (with $b=\frac{3}{8}$) and the
Sobolev embedding $H^{\frac{3}{2}}(M)\subset L^{\infty}(M)$, we obtain a second bound for $I(L,M)$
\begin{eqnarray*}
I(L,N) 
& \leq & 
\prod_{j=0}^{1}
\|u_{j}^{L_j N_j}\|_{L^{4}(\R\,;\,L^{2}(M))}
\,
\prod_{j=2}^{3}
\|u_{j}^{L_j N_j}\|_{L^{4}(\R\,;\,L^{\infty}(M))}
\\
& \lesssim &
\frac{(L_0 L_1 L_2 L_3)^{\frac{3}{8}}}
{L_0^{b'}(L_1 L_2 L_3)^{b}}
\frac{N_0^s}{(N_1 N_2 N_3)^{s}}
(N_2 N_3)^{\frac{3}{2}}
\prod_{j=0}^{3}c_{j}(L_j,N_j).
\end{eqnarray*}
We now state a technical lemma, which is useful for the interpolation argument.
\begin{lemme}\label{le2.3}
For every $s>s_{0}$, there exists $(b,b')$ satisfying (\ref{star}), $\epsilon >0$, 
and $\theta\in]0,1[$ such that 
\begin{equation}\label{eq2.12}
s>\frac{3}{2}\theta+s_{0}(1-\theta),
\end{equation}
\begin{equation}\label{eq2.13}
b'>\frac{3}{8}\,\theta+\big(\frac{1}{2}+\epsilon\big)(1-\theta).
\end{equation}

\end{lemme}
{\bf Proof.}
For a fixed $s>s_{0}$, we shall choose the parameters following the scheme
$$
\theta \rightarrow \epsilon \rightarrow (b,b').
$$
We first choose $\theta \in]0,1[$ such that \eqref{eq2.12} holds: 
$$\theta (\frac 3 2 - s_{0}) < s-s_{0}$$
(any choice works if $s_{0}\geq \frac 3 2$, otherwise $\theta$ has to be close enough to $0$). Next, we choose $\epsilon>0$ such that
\begin{equation}\label{eq2.15}
\epsilon<\frac{\theta}{8(3-\theta)}.
\end{equation}
We finally set
\begin{equation}\label{star'}
b:=\frac{1}{2}+\epsilon,\quad b':=\frac{1}{2}-2\epsilon.
\end{equation}
Clearly (\ref{star'}) ensures (\ref{star}) for $\epsilon$ small enough.
It remains finally to observe that (\ref{eq2.13}) follows from (\ref{eq2.15}).
\\

Using the technical lemma and interpolating between the two bounds for $I(L,N)$ with respective weights $1-\theta$ and $\theta$,
we obtain that there exist $\gamma_1>0$ and $\gamma_2>0$ such that
$$
I(L,N)\lesssim
\left(\frac{N_0}{N_1}\right)^{s}\frac{1}{(N_2 N_3)^{\gamma_1} (L_0 L_1 L_ 2 L_ 3)^{\gamma_2}}
\prod_{j=0}^{3}c_{j}(L_j,N_j).
$$
The summation over $L_0$, $L_1$, $L_2$, $L_3$, $N_2$, $N_3$ can be performed via a crude argument 
of summation of geometric series which gives
$$
\sum_{L}\sum_{N}
I(L,N)\lesssim
\|w_2\|_{L^2}
\|w_3\|_{L^2}\,
\sum_{N_0\leq 3C N_1}
\left(\frac{N_0}{N_1}\right)^{s}\alpha(N_1)\,\beta(N_0),
$$
where
$$
(\alpha(N_1))^{2}=\sum_{L_1}(c_{1}(L_1, N_1))^{2},\quad (\beta(N_0))^{2}=\sum_{L_0}(c_0(L_0, N_0))^{2}.
$$
Suppose now that $N_{1}=2^{l}N_0$, where $l\geq -l_0$ where $l_0$ is a fixed positive integer depending only on $3C$.
We then write, via the Cauchy-Schwarz inequality in $N_0$,
\begin{eqnarray*}
\sum_{N_0\leq 3C N_1}
\left(\frac{N_0}{N_1}\right)^{s}
\alpha(N_1)\beta(N_0)
& = &
\sum_{l\geq -l_0}\,
\sum_{N_0}
2^{-sl}
\,
\alpha(2^{l}N_0)\,\beta(N_0)
\\
& \lesssim &
\sum_{l\geq -l_0}
2^{-sl}
\Big\{
\sum_{N_0}
(\alpha(2^{l}N_0))^{2}
\Big\}^{\frac{1}{2}}
\Big\{
\sum_{N_0}
(\beta(N_0))^{2}
\Big\}^{\frac{1}{2}}
\\
& \lesssim &
\|w_1\|_{L^2}\|w_2\|_{L^2}.
\end{eqnarray*}
This completes the proof of Proposition \ref{prop2.2}.
\qed
\subsection{The contraction argument}
We are finally in position to prove Theorem~\ref{th3}. Denote by $S(t) = e^{it\Delta}$ the free evolution. The function $u\in C^{\infty}(\R\times M)$ solves (\ref{nls}) if and only if it also solves
the integral equation (Duhamel form)
\begin{equation}\label{Duhamel}
u(t)=S(t)u_{0}- i\int_{0}^{t}S(t-\tau)\{|u(\tau)|^{2}u(\tau)\}d\tau.
\end{equation}
The next proposition contains the basic linear estimates.
\begin{proposition}\label{linear-estimates}
Let $\psi\in C_{0}^{\infty}(\R)$ be equal to $1$ on $[-1,1]$. Then
\begin{equation}\label{eq2.16}
\|\psi(t)S(t)u_{0}\|_{X^{s,b}(\R\times M)}\lesssim\|u_0\|_{H^s(M)}
\end{equation}
and 
\begin{equation}\label{eq2.17}
\|\psi(t/T)\int_{0}^{t}S(t-\tau)F(\tau)d\tau\|_{X^{s,b}(\R_t\times M)}\lesssim T^{1-b-b'}\|F\|_{X^{s,-b'}(\R\times M)},
\end{equation}
provided 
\begin{eqnarray}\label{restrictions}
0<b'<\frac{1}{2}<b,\quad b+b'<1,\quad 0< T \leq 1.
\end{eqnarray}
\end{proposition}
\proof
Recall that
\begin{eqnarray}\label{describtion}
\|u\|_{X^{s,b}(\R\times M)}=\|S(-t)u(t)\|_{H^{b}_{t}(\R\,;\,H^{s}(M))}
\end{eqnarray}
which yields the equality
$$
\|\psi(t)S(t)u_{0}\|_{X^{s,b}(\R\times M)}=\|\psi\|_{H^{b}(\R)}\|u_0\|_{H^s(M)}
$$
and the proof of (\ref{eq2.16}) is therefore completed.
Next we prove (\ref{eq2.17}). Under the assumption (\ref{restrictions}), the elementary inequality~\cite[(3.11)]{Gi95} reads 
\begin{eqnarray}\label{basic}
\|\psi(t/T)\int_{0}^{t}g(\tau)d\tau\|_{H^{b}(\R)}\lesssim 
T^{1-b-b'}\|g\|_{H^{-b'}(\R)}.
\end{eqnarray}
We now show that (\ref{basic}) implies  (\ref{eq2.17}). First, due to (\ref{describtion}), we observe that 
(\ref{eq2.17}) is equivalent to
\begin{equation}\label{novo}
\|\psi(t/T)\int_{0}^{t}G(\tau)d\tau\|_{H^{b}(\R_t\,;\,H^{s}(M))}
\lesssim T^{1-b-b'}\|G\|_{H^{-b'}(\R\,;\,H^{s}(M))}.
\end{equation}
Next, using (\ref{basic}), we obtain that
\begin{equation}\label{pointwise}
\|\psi(t/T)\int_{0}^{t}P_{k}\,G(\tau)d\tau\|_{H^{b}(\R_t)}
\lesssim T^{1-b-b'}\|P_k\, G\|_{H^{-b'}(\R_t)},
\end{equation}
pointwise on $M$, if $G\in C_{0}^{\infty}(\R\times M)$.
Squaring (\ref{pointwise}), integration over $M$, multiplying it with $\langle\mu_k\rangle^{s}$ and finally summing over $k$ gives 
(\ref{novo}).
\qed
\\
\newline
Next, for $T>0$, we define the restriction space 
$X^{s,b}_{T}:=X^{s,b}([-T,T]\times M)$, equipped with the norm
$$
\|u\|_{X^{s,b}_{T}}=\inf_{w\in X^{s,b}}
\{
\|w\|_{X^{s,b}},\quad{\rm with }\quad w|_{[-T,T]}=u
\}.
$$
Note that if $b>1/2$, the space $X^{s,b}_{T}$ is continuously embedded in $C([-T,T]\,;\,H^{s}(M))$. 
Using (\ref{eq2.16}), we get
\begin{equation}\label{eq2.16-bis}
\|\psi(t)S(t)u_{0}\|_{X^{s,b}_{T}}\lesssim\|u_0\|_{H^s},
\end{equation}
if $T\leq 1$. Moreover propositions \ref{linear-estimates} and \ref{prop2.2} yield
\begin{eqnarray}\label{jerry}
\quad
\|\int_{0}^{t}S(t-\tau)u_1(\tau)\, \overline{u_2}(\tau)\, u_3(\tau)d\tau\|_{X^{s,b}_{T}}\lesssim 
T^{1-b-b'}\|u_1\|_{X^{s,b}_T}\|u_2\|_{X^{s,b}_T}\|u_3\|_{X^{s,b}_T}\, ,
\end{eqnarray}
provided that (\ref{restrictions}) holds.
Let now $u(t)$ and $v(t)$ be two smooth solutions of the cubic Schr\"odinger equation such that $u(0)$ and $v(0)$ belong
to a fixed bounded set $B$ in $H^{s}(M)$. Using \eqref{Duhamel}, (\ref{eq2.16-bis}) and (\ref{jerry}), we obtain that there exist
$T$ and $C>0$, depending only on $B$, such that
\begin{eqnarray}\label{bound}
\|u\|_{X^{s,b}_{T}}+\|v\|_{X^{s,b}_{T}}\leq C.
\end{eqnarray}
Writing $|u|^{2}u-|v|^{2}v=u^{2}(\overline{u}-\overline{v})+\overline{v}(u-v)(u+v)$ using again 
(\ref{eq2.16-bis}) and (\ref{jerry}) and the bound (\ref{bound}), we obtain that
$$
\|u-v\|_{X^{s,b}_T}\lesssim C\|u(0)-v(0)\|_{H^s}.
$$
Since $b>1/2$, we finally get
$$
\|u-v\|_{C([-T,T]\,;\,H^s(M))}\lesssim C\|u(0)-v(0)\|_{H^s}.
$$
which yields the uniform wellposedness.
\qed
\begin{remarque}
It is worth noticing that all the results we know about uniform wellposedness for (\ref{nls}) fit in theorem \ref{th3}
in the sense that (\ref{eq.bilin}) is satisfied. For instance, by combining 
inequality (\ref{strichartz}) for $p$ close to $1/2$ with Sobolev embedding one obtains easily
(\ref{eq.bilin}) for every $s_0>1/2$, which yields the uniform wellposedness result of
\cite{BuGeTz01-1}.
On the other hand, we do not know whether uniform wellposedness in $H^s$ implies (\ref{eq.bilin}) with $s_0=s$.
However, we present a simple argument below which shows that, 
if the flow map (\ref{flow}) is $C^3$ near the origin, then  (\ref{eq.bilin})
holds with $s_0=s$. 
Notice that conversely, the smoothness of the flow map is also a consequence of the arguments of this section.
\\
Assume the flow map is $C^3$. Then its third differential at the origin is the trilinear operator
\begin{multline*}
{\cal T}(h_1, h_2,h_3)=
-2i\int_{0}^{t}S(t-\tau)\big[
S(\tau)h_1\, \overline{S(\tau)h_2}\,S(\tau)h_3
+\\
S(\tau)h_2\, \overline{S(\tau)h_3}\,S(\tau)h_1
+
S(\tau)h_3\, \overline{S(\tau)h_1}\,S(\tau)h_2
\big]
d\tau.
\end{multline*}
For a fixed $t=T$, we take the $H^s$ scalar product of ${\cal T}(f,f,g)$ with $g$. Then the continuity of ${\cal T}$ implies
$$
\left|\int_{0}^{T}\int_{M}
\big[\,
2\,|S(\tau)f\, S(\tau) g|^2+(S(\tau)f)^{2}\,(\overline{S(\tau)g})^{2}
\big]
d\tau
\right|
\leq
C\|f\|_{H^s}^{2}
\|g\|_{H^s}\|g\|_{H^{-s}}
$$
Taking $f$ and $g$ satisfying (\ref{lokalizatzia}) with $N\leq L$, we obtain (\ref{eq.bilin}) with $s_0=s$.
\end{remarque} 
%%%%%%%%%%%%%%%%%%%%%%%%%%%%%%%%%%%%%%%%%%%%%%%%%%%%%%%%%%%%%%%%%%%%%%%%%%%%%%%%%%%%%%%%%%%%%%%%%%%%%%%%%%%%%%%%%%%%%%%%%%%%%%%%%%%%%
\section{From spectral estimates to evolution estimates}
The main result in this section is the proof of Theorem~\ref{th1} (assuming Theorem~\ref{th2}). Taking benefit of Theorem~\ref{th3}, 
it suffices to prove the bilinear Strichartz estimate~\eqref{eq.bilin} for $s_{0}>1/4$. Using Lemma~\ref{le2.5}, this in turn 
reduces to the proof of~\eqref{eq2.2}. On the sphere $\mathbb{S}^2$ the situation is a little simpler since we know exactly the 
spectrum of the Laplace operator ($\mu_{k}= k(k+1)$, $k\in\N$). We can use this knowledge to show that the 
bilinear estimate~\eqref{eq2.2} for $s_{0}>1/4$ can be deduced from  Theorem~\ref{th2}. Then we  extend this result to the case of 
Zoll manifolds. In this case the localization of the spectrum is still well understood (see Proposition~\ref{prop3.3}) and, by introducing 
a suitable abstract perturbation of the Laplace operator, we are able to reduce the study to the case of the sphere. 
In the case of the torus, we did not find inequality \eqref{eq.bilin} explicitly written in the literature.
Nevertheless, it can be easily deduced (for any $s_{0}>0$) from the analysis in \cite{Bo93} as we show at the end of that section. 
%%%%%%%%%%%%%%%%%%%%%%%%%%%%%%%%%%%%%%%%%%%%%%%%%%%%%%%%%%%%%%%%%%%%%%%%%%%%%%%%%%%%%%%%%%%%%%%%%%%%%%%%%%%%%%%%%%%%%%%%%%5
\subsection{The case of the sphere $\mathbb{S}^2$}
We are now going to prove that the bilinear estimate on the two dimensional sphere endowed with its standard metric is a 
consequence of the precise knowledge of the spectrum: $\mu_{k}= k(k+1)$, $k\in\N$ and Theorem~\ref{th2}:
\begin{proposition}\label{prop3.1}
Let $s>1/4$. There exists $C>0$ such that for any $u_0, v_0 \in L^2(\mathbb{S}^2 )$ satisfying
\begin{eqnarray}\label{hipoteza}
\11_{\sqrt{-\Delta}\in [N,2N]}(u_0)=u_0,\quad
\11_{\sqrt{-\Delta}\in [L,2L]}(v_0)=v_0
\end{eqnarray}
one has 
\begin{equation}\label{sfera}
\|e^{it \Delta}u_0\,\, e^{it\Delta} v_0\|_{L^2((0,1)_t\times \mathbb{S}^2)}
\leq 
C(\min(N,L))^{s}\|u_{0}\|_{L^2(\mathbb{S}^2)} \|v_{0}\|_{L^2(\mathbb{S}^2)}.
\end{equation}
\end{proposition}
\proof
We write
$$
e^{it \Delta}u_0=\sum_{k=0}^{\infty}e^{-itk(k+1)}P_{k}(u_0)\quad
e^{it \Delta}u_0=\sum_{k=0}^{\infty}e^{-itk(k+1)}P_{k}(v_0),
$$
where $P_k$ is the projector on spherical harmonics of degree $k$.
Using (\ref{hipoteza}), we get
$$
e^{it \Delta}u_0\,\, e^{it\Delta} v_0
=\sum_{N\leq k\leq 2N}\,\,\sum_{L\leq l\leq 2L}e^{-it(k(k+1)+l(l+1))}P_{k}(u_0)P_{l}(v_0).
$$
Therefore using Parseval identity, 
\begin{eqnarray*}
\Big\|
e^{it \Delta}u_0\,\, e^{it\Delta} v_0
\Big\|^{2}_{L^{2}([0,2\pi]\times \mathbb{S}^2)} & = & \sum_{\tau=0}^{\infty}
\Big\|\sum_{\tau=k(k+1)+l(l+1)}P_{k}(u_0)P_{l}(v_0)\Big\|^{2}_{L^{2}(\mathbb{S}^2)}
\\
& \leq &
\sum_{\tau=0}^{\infty}\alpha_{NL}(\tau)\sum_{\tau=k(k+1)+l(l+1)}
\|P_{k}(u_0)P_{l}(v_0)\|^{2}_{L^{2}(\mathbb{S}^2)},
\end{eqnarray*}
where
$$
\alpha_{NL}(\tau)=\# \,\,\{ (k,l)\in \N\times\N,\quad N\leq k\leq 2N,\,\, L\leq l\leq 2L,\,\,k(k+1)+l(l+1)=\tau\}.
$$
In order to bound $\alpha_{NL}(\tau)$, we need a lemma giving a bound on the number of lattice points
on an arc of a circle in $\R^2$.
\begin{lemme}\label{le3.2}
Let $M$ and $N$ be two positive integers.
Then for any $\epsilon>0$ there exists $C>0$ such that
$$
\#\{(k_1,k_2)\in\mathbb{N}^2\,:\, N\leq k_{1} \leq 2N,\quad k_1^2+k_2^2=M\}
\leq C N^{\epsilon}.
$$\end{lemme}
\proof
We consider two cases.
If $M\leq N^{100}$ then we use the elementary bound on the number of divisors of $M$ in the ring of Gauss integers,
\begin{equation}\label{gaussdiv}
\# \{ z\in \mathbb{Z} [i]\,:\,   \exists \widetilde z \in \mathbb{Z}[i], z\widetilde z = M\} =\mathcal{O} (M^\delta), \forall \delta >0,
\end{equation}
which gives a bound $\mathcal {O} (N^ \epsilon)$.
If $M>N^{100}$ then
we simply observe that $k_2$ takes values in the interval
$[\sqrt{M-4N^2},\sqrt{M-N^2}]$ which in the considered case
contains at most one integer.
The proof of the lemma is completed taking into account that if we fix
$k_2$ then $k_1$ can not take more than one value.
\qed 
\\
\newline
Noticing that 
\begin{eqnarray}\label{diofant}
k(k+1)+ l(l+1)= \tau\Leftrightarrow 4\tau+2= (2k+1)^2 + (2l+1) ^2,
\end{eqnarray}
and using Lemma~\ref{le3.2}, we obtain the bound
$$
\sup_{\tau}\alpha_{NL}(\tau)\leq C_{\epsilon}[\min(N,L)]^{\epsilon}
$$
for any $\epsilon>0$. Therefore using Theorem \ref{th2}, we obtain that
$$
\Big\|
e^{it \Delta}u_0\,\, e^{it\Delta} v_0
\Big\|_{L^{2}([0,2\pi]\times \mathbb{S}^2)}
\leq C_{\epsilon}[\min(N,L)]^{\frac{1}{4}+\epsilon}
\|u_0\|_{L^2(\mathbb{S}^2)}\|v_0\|_{L^2(\mathbb{S}^2)}.
$$
Inequality (\ref{sfera}) now follows simply by restricting the left hand-side in the last estimate to $(0,1)_t$
and taking $\epsilon:=s-1/4$.
\qed
%%%%%%%%%%%%%%%%%%%%%%%%%%%%%%%%%%%%%%%%%%%%%%%%%%%%%%%%%%%%%%%%%%%%%%%%%%%%%%%%%%%%%%%%%%%%%%%%%%%%%%%%%%%%%%%%%%%%%%%
%%%%%%%%%%%%%%%%%%%%%%%%%%%%%%%%%%%%%%%%%%%%%%%%%%%%%%%%%%%%%%%%%%%%%%%%%%%%%%%%%%%%%%%%%%%%%%%%%%%%%%%%%%%%%%%%%%%%%%%
\subsection{The case of Zoll surfaces}
In this section we extend Proposition~\ref{prop3.1} to the case of Zoll surfaces. A Zoll manifold is a Riemannian manifold such that 
all geodesics are closed with a common period. The starting point is the following result 
(see Guillemin~\cite{Gu77} and Colin de Verdi\`ere~\cite{CdV79}).
\begin{proposition}\label{prop3.3} 
If the geodesics of $M$ are $2\pi$ periodic, 
there exists $\alpha \in \mathbb{N}$ and $E>0$ such that the spectrum of $-\Delta$ is contained in $\cup_{k=1}^{+\infty} I_k$ where
$$I_k= \left[ \left(k+\frac \alpha 4 \right)^2-E,\, \left( k + \frac \alpha 4 \right)^2 + E\right].$$
\end{proposition}
In fact we are going to reduce the analysis here  to the analysis already performed in the previous section.
Indeed, denote by $(e_n, \lambda_n)$ the sequence of eigenfunctions of $-\Delta$ associated to the eigenvalues 
$\lambda_n$, fix $N_0$ such that for $n\geq N_0$ the eigenvalue $\lambda _n$ is exactly in one interval $I_k$ 
(it is possible since for large $k$ these intervals are disjoint). Define the abstract operator $\widetilde{\Delta}$ on $L^2(M)$ 
by the relations:
\begin{equation}
-\widetilde {\Delta} e_n 
= 
\begin{cases}
\lambda_n e_n &\text{ if } n\leq N_0
\\
\left(k+\frac \alpha 4 \right)^2 e_n &\text{ if } n>N_0 \text{ and } \lambda_n \in I_k
\end{cases}
\end{equation}
To show that Proposition~\ref{prop3.1} holds with $\mathbb{S}^2$ replaced by $M$, it is enough, according to Lemma~\ref{le2.5},  
to prove~\eqref{eq2.2} for $s_0= 1/4 +\varepsilon$ and $X^{0,b}_{-\Delta}$.
But according to Remark~\ref{re2.1}, 
$$
X^{0,b}_{-\Delta}(\R\times M)=X^{0,b}_{-\widetilde {\Delta}}(\R\times M).
$$
Using Lemma~\ref{le2.5} in the other way, it is enough to show that~\eqref{eq2.1} is fulfilled for $s_0= 1/4+ \epsilon$ 
and $A= - \widetilde{\Delta}$. Taking into account that 
\begin{equation}\label{eq3.4}
\11_{\left[\left(k+\frac \alpha 4 \right)^2-E, \left(k+\frac \alpha 4 \right)^2+E \right]}( - \Delta)
=
\11_{\left[\left(k+\frac \alpha 4 \right)^2-E, \left(k+\frac \alpha 4 \right)^2+E \right]}( - \widetilde{\Delta})
\end{equation}
we see that for these spectral projectors, the estimates in Theorem~\ref{th2} hold. 
\\
Replacing, in the proof of Proposition~\ref{prop3.1}, projectors $P_k$ by (\ref{eq3.4}) and 
(\ref{diofant}) by
\begin{eqnarray*}
\left(k+\frac \alpha 4 \right)^2+\left(l+\frac \alpha 4 \right)^2=\tau
\Leftrightarrow 
(4k+\alpha)^{2}+(4l+\alpha)^{2}=16\tau,
\end{eqnarray*}
we see that the proof holds equally well for $-\widetilde {\Delta}$.
\qed 
%%%%%%%%%%%%%%%%%%%%%%%%%%%%%%%%%%%%%%%%%%%%%%%%%%%%%%%%%%%%%%%%%%%%%%%%%%%%%%%%%%%%%%%%%%%%%%%%%%%%%%%%%%%%%%%%%%%%%%%%
\subsection{The case of the torus $\mathbb{T}^2$}
In this section $\mathbb{T}^n= \mathbb{R}^n / \mathbb{Z}^n $. 
We first recall  the following lemma due to Bourgain.
\begin{lemme} Let $Q_0=[-1,1]^2$ be the unit square in $\mathbb{R}^2$ and 
$Q= (a,b)+ A Q_0$, $a,b\in \mathbb{Z}$ be a square in $\mathbb{R}^2$ of size $A\in\N$ .
Then for any $\epsilon>0$ there exists $C_{\epsilon}>0$ such that for any
$$
u(t,x_1,x_2)= \sum_{(n_1, n_2)\in Q\cap\Z^2} c_{n_1, n_2}\,\, e^{2\pi i(t (n_1^2+ n_2^2)+ x_1 n_1 +  x_2 n_2)}
$$
one has 
\begin{equation}
\|u\|_{L^4(\mathbb{T}^{3})}\leq \,
C_{\epsilon}\,A^\epsilon 
\Big( \sum_{n_1, n_2} |c_{n_1, n_2}|^2\Big)^{1/2}
\end{equation}
\end{lemme}
\proof
For the sake of completeness, we give the proof. Write $(n_1, n_2) = (a,b) +(q_1, q_2)$. Then
$$
u(t,x_1,x_2)=  e^{2\pi it(a^2+b^2)}
\sum_{(q_1, q_2)\in AQ_0\cap \Z^2} c_{q_1, q_2}\,\, e^{2\pi i(t (q_1^2+ q_2^2)+ (x_1+2ta) q_1 +  (x_2+2tb) q_2)}
$$
and the change of variables $(x_1,x_2,t)\mapsto (x_1-2ta,x_2-2tb,t)$ shows that it is enough to prove the result for $a=b=0$. Then
$$u^2(t,x_1,x_2)= 
\sum_{|q_1|, |q_2|, |r_1|, |r_2|\leq A} c_{q_1,q_2}\, c_{r_1,r_2} 
\,\,
e^{2\pi it(q_1^2+ q_2^2+ r_1^2+ r_2^2)}e^{2\pi ix_1(q_1+r_1)}e^{2\pi ix_2(q_2+r_2)}$$
and with 
$$
I_{\tau,p_1,p_2}:=\{(q_1,q_2,r_1,r_2)\,:\,  q_1^2+ q_2^2+ r_1^2+r_2^2= \tau,\, 
q_j+r_j=p_j,\, |q_j|\leq A,\, |r_j|\leq A,\, j=1,2\},
$$
using the Parseval identity and the Cauchy-Schwarz inequality, we get
\begin{multline}
\begin{aligned}
\|u^2\|_{L^2(\T^3)} ^2 &= 
\sum_{\tau, p_1, p_2} \Bigl| 
\sum_{(q_1,q_2, r_1, r_2)\in I_{\tau, p_1, p_2}}
c_{q_1,q_2}
\,
c_{r_1, r_2}
\Bigr|^{2}
\\ 
 &\leq
\sup_{\tau,p_1,p_2}
|I_{\tau, p_1, p_2}|
\Big(
\sum_{q_1,q_2}|c_{q_1,q_2}|^2
\Big)
\Big(
\sum_{r_1,r_2}|c_{r_1,r_2}|^2
\Big),
\end{aligned}
\end{multline}
where $|I_{\tau, p_1, p_2}|$ denotes the number of elements of $I_{\tau, p_1, p_2}$.

We now estimate $|I_{\tau, p_1, p_2}|$. Fix $(\tau, p_1, p_2)\in\Z^3$. Let $(q_1,q_2,r_1,r_2)\in I_{\tau, p_1, p_2}$.
Clearly $|\tau|\leq 4A^2$ and $|p_{j}|\leq 2A$, $j=1,2$. Observing that
$$\tau = q_1^2+q_2^2+(p_1-q_1)^2+ (p_2-q_2)^2\Leftrightarrow 
2\tau- p_1^2-p_2^2= (2q_1-p_1)^2+ (2q_2-p_2)^2$$
and using the divisor bound in the ring $\Z[i]$ of Gauss integers~\eqref{gaussdiv}, 
we infer that for any $\epsilon>0$ there exists $C_{\epsilon}$ such that
$|I_{\tau, p_1, p_2}|\leq C_{\epsilon}A^{\epsilon}$, which completes the proof.
\qed
\\

Following the ideas of Bourgain, we now deduce the following bilinear estimate,
which readily implies (\ref{eq.bilin}) in the case of torus with $s_0=\epsilon$, $\epsilon$ being an arbitrary positive number.  
\begin{proposition}\label{carre}
For any $\epsilon>0$ there exists $C_{\epsilon}>0$ such that for every couple $N,L$ of dyadic integers one has
$$
\|e^{it \Delta}u_0\,\, e^{it\Delta} v_0\|_{L^2( \mathbb{T}^3)}\leq 
C_{\epsilon}\, (\min(N,L))^\epsilon\,\|u_0\|_{L^2( \mathbb{T}^2)}\|v_0\|_{L^2( \mathbb{T}^2)}\, ,
$$
if
$$
u_0(x_1,x_2)= \sum_{N\leq \sup (|n_1|,|n_2|)< 2N}c_{n_1,n_2}e^{2\pi i(x_1n_1+x_2n_2)},
$$
and
$$
v_0(x_1,x_2)= \sum_{L\leq\sup( |m_1|,|m_2|)< 2L}d_{m_1,m_2}e^{2\pi i(x_1m_1+x_2m_2)}\,\,.
$$
\end{proposition}
\proof
We can suppose $N\leq L/10$, the case $L\leq N/10$ being similar and the case $N\sim L$ being a 
consequence of H\"older inequality and the previous lemma. Then we can cover the shell
$$
\{(m_1, m_2)\in\Z^2\,:\, L\leq \sup(|m_1|, |m_2|)< 2L\}
$$ 
by squares $Q_\alpha$ with disjoint interiors  :
$$
Q_\alpha=\{(m_1, m_2)\in\Z^2\,:\, \sup( |m_1- \alpha_1|, |m_2- \alpha_2|)\leq N\},
$$
where $\alpha= (\alpha_1,\alpha_2)$ belongs to a suitable subset of $\Z^2$.
Correspondingly we associate to this decomposition the natural splitting
\begin{eqnarray*}
v(t,x_1,x_2)
& = & \sum_\alpha \sum_{(m_1, m_2)\in Q_\alpha}\,d_{m_1,m_2}\,e^{2\pi i(t(m_1^2+ m_2^2) + x_1m_1+x_2m_2)}\,
\\
& := &
\sum_{\alpha} v_{\alpha}(t,x_1,x_2)
\end{eqnarray*}
and we have 
\begin{equation}\begin{aligned}
\|u\,v\|_{L^2( \mathbb{T}^3)}^2 
&= \|\sum _\alpha u v_\alpha\|_{L^2( \mathbb{T}^3)}^2 \leq C \sum_\alpha \|uv_\alpha\|_{L^2( \mathbb{T}^3)}^2\\
&\leq C \sum_\alpha \|u\|_{L^4( \mathbb{T}^3)}^2\|v_\alpha\|_{L^4( \mathbb{T}^3)}^2 \leq C \sum_\alpha N^\epsilon 
\|u_0\|_{L^2( \mathbb{T}^2)}^2 \|v_{\alpha}(0)\|_{L^2( \mathbb{T}^2)}^2\\
&\leq C N^\epsilon \|u_0\|_{L^2( \mathbb{T}^2)}^2 \|v_0\|_{L^2( \mathbb{T}^2)}^2
\end{aligned}
\end{equation}
where we have used for the first inequality that the Fourier series expansion of 
$uv_\alpha$ is localized in the square of center $\alpha$ and size $2N$ and consequently the terms in the series
$\sum_\alpha uv_\alpha$ are quasi orthogonal in $L^2$.
\qed  
%%%%%%%%%%%%%%%%%%%%%%%%%%%%%%%%%%%%%%%%%%%%%%%%%%%%%%%%%%%%%%%%%%%%%%%%%%%%%%%%%%%%%%%%%%%%%%%%%%%%%%%%%%%%%%%%%%%%%%%%%%%%%%% 
\section{Bilinear eigenfunction estimates}\label{sec4}
In this section $M$ is any compact Riemannian surface and $\Delta$ the Laplace operator on functions on $M$. 
We are going to prove Theorem~\ref{th2}. The strategy consists in performing suitable canonical transformations to reduce 
the study to an oscillatory integral which can be studied by a bilinearization of an argument by H\"ormander~\cite{Ho73}. 
\subsection{Reduction to an oscillatory integral}
To prove Theorem~\ref{th2}, it is enough to consider the case when
$\chi\in{\cal S}(\R)$ and its Fourier transform is supported in $[\varepsilon,2\varepsilon]$
 with $\varepsilon>0$ small. The starting point is Sogge's Proposition~\ref{propsogge}, which gives
$$
\chi (\sqrt{-\Delta}-\lambda)f
=
\lambda^{\frac{1}{2}}{\mathcal {T}}_{\lambda}f+R_{\lambda}f,
$$
with
$$
\|R_{\lambda}f\|_{L^{\infty}}\leq C\|f\|_{L^2}
$$
and, in a coordinate system close to $x_0\in M$,
$$
{\mathcal {T}}_{\lambda}f(x)=\int_{\R^2}e^{i\lambda\varphi(x,y)}a(x,y)f(y)dy
$$
with $a(x,y)\in C_{0}^{\infty}$ supported in
$$
\frac{\varepsilon}{C}
\leq |x-y|
\leq 
C\varepsilon,
$$
and $\varphi(x,y)=-d_{g}(x,y)$. Notice that here we kept only the main term of the symbol $a(x,y,\lambda)$.\par
Consequently it is enough to prove
$$
\|{\mathcal {T}}_{\lambda}f\,\,\mathcal{T}_{\mu}g\|_{L^2}
\leq
\frac{C}{\lambda^{\frac{1}{4}}\,\mu^{\frac{1}{2}}}
\|f\|_{L^2}\|g\|_{L^2}
$$
 if $1\leq\lambda\leq\mu$. Using partitions of unity, we can suppose that $x_{0}=0$ and $a(x,y)$
is supported in the set $|x|\leq\delta\ll\frac{\varepsilon}{C}$.
We introduce geodesic polar coordinates centered at $0$ by setting $y=\exp_{o}(re^{it})$, so that
$$
{\mathcal {T}}_{\lambda}f(x)=
\int_{0}^{\infty}\int_{0}^{2\pi}
e^{i\lambda\varphi_{r}(x,t)}
a_{r}(x,t)f_{r}(t)drdt
$$
with $\varphi_{r}(x,t)=\varphi(x,y)$, $a_{r}(x,t)=a(x,y)$ and $f_{r}(t)=f(y)$.
Denote by
$$
{\mathcal {T}}_{\lambda}^{r}F(x)=
\int_{0}^{2\pi}
e^{i\lambda\varphi_{r}(x,t)}
a_{r}(x,t)F(t)dt
$$
which gives
$$
{\mathcal {T}}_{\lambda}f(x)\,\,{\mathcal{T}}_{\mu}g(x)
=
\int_{\varepsilon/C}^{C\varepsilon}
\int_{\varepsilon/C}^{C\varepsilon}
({\mathcal {T}}_{\lambda}^{r}f_{r})(x)\,
({\mathcal{T}}_{\mu}^{q}g_{q})(x)drdq
$$
and hence
$$
\|{\mathcal {T}}_{\lambda}f\,{\mathcal{T}}_{\mu}g\|_{L^2}
\leq
\int_{\varepsilon/C}^{C\varepsilon}
\int_{\varepsilon/C}^{C\varepsilon}
\|
{\mathcal {T}}_{\lambda}^{r}f_r\,\,{\mathcal{T}}_{\mu}^{q}g_q
\|_{L^2}
drdq\, .
$$
Consequently we have reduced the analysis to proving
\begin{equation}\label{eq4.1}
\|{\mathcal {T}}_{\lambda}^{r}F\,{\mathcal{T}}_{\mu}^{q}G\|_{L^2}
\leq
\frac{C}{\lambda^{\frac{1}{4}}\,\mu^{\frac{1}{2}}}
\|F\|_{L^2([0,2\pi])}\|G\|_{L^2([0,2\pi])}
\end{equation}
uniformly for $1\leq\lambda\leq \mu$,
$\frac{\varepsilon}{C}\leq r,q\leq C\varepsilon$.
\\
\newline 
It turns out that $\varphi_{r}(x,t)$ is a Carleson-Sj\"olin phase (see \cite{CaSj72, Ho73}).
\begin{lemme}
\label{le.carsjhyp}
For every $0<\delta<1$ there exist $c>0$, $\varepsilon>0$ such that for every $|x|\leq \varepsilon$,
$$
\left|\det
[\nabla _{x}\partial _{t} \varphi_{r}(x,t), \nabla_{x}\partial _{s}\varphi_{q}(x, s)]
\right| 
\geq c\, ,
$$
if
$|t-s|\geq \delta$ and $|t+\pi-s|\geq \delta$. In addition, for every $t\in[0,2\pi]$,
$$
\left|\det
[\nabla _{x}\partial _{t} \varphi_{r}(x,t), \nabla_{x}\partial ^2_{t}\varphi_{r}(x, t)
]\right| \geq c\, .
$$
\end{lemme} 
\proof
Let $u=u(x,y)\in T_{y}M$ be the unit vector such that
$$
\exp_{y}(-\varphi(x,y)u(x,y))=x.
$$
Differentiating with respect to $x$ this identity, we get for $x=0$, and any $h\in T_{o}M$,
\begin{equation}\label{eq4.100}
h=
T_{ru(0,y)}(\exp_{y})
\left[
-d_{x}\varphi(0,y)\cdot hu(0,y)-r\, T_{x}u(0,y)\cdot h
\right].
\end{equation}
Remark that since $y=\exp_{o}(re^{it})$,
\begin{equation}
T_{ru(0,y)}(\exp_{y})\cdot u(0,y)=-e^{it}
\end{equation}
Take the scalar product with $e^{it}$ in~\eqref{eq4.100}. Using Gauss' Lemma (see for example~\cite[3.70]{GaHuLa90}), we get
$$
d_{x}\varphi(0,y)\cdot h=g_{o}(h,e^{it}),
$$
i.e.
$$
\nabla_{x}\varphi(0,y)=e^{it}.
$$
Differentiating with respect to $t$, this implies
$$
\nabla_{x}\partial_{t}\varphi_{r}(0,t)=ie^{it}
$$
so that
$$\left|\det
\left[
\nabla_{x}\partial_{t}\,\varphi_{r}(0,t),
\nabla_{x}\partial_{s}\varphi_{q}(0,s)
\right]\right|= \left|\sin(t-s)\right| \neq 0,$$
$$
\det
\left[
\nabla_{x}\partial_{t}\,\varphi_{r}(0,t),
\nabla_{x}\partial^{2}_{t}\varphi_{r}(0,t)
\right]
\neq 0,
$$
for any $t,s $ satisfying the assumptions of the lemma. By continuity this remains true for $x$ close enough to $0$. 
\endproof
Let us come back to the proof of inequality (\ref{eq4.1}). Write
$$
{\mathcal {T}}_{\lambda}^{r}f(x)\,{\mathcal{T}}_{\mu}^{q}g(x)
=
\int\int_{[0,2\pi]^2}
e^{i(\lambda\varphi_{r}(x,t)+\mu\varphi_{q}(x,s))}
a_{r}(x,t){a}_{q}(x,s)f(t)g(s)dtds.
$$
First, let us assume that, if $(t,s)\in {\rm supp }(f)\times {\rm supp }(g)$, then $(t,s)$ satisfies
$|t-s|\geq \delta$ and $|t+\pi-s|\geq \delta$.
We may moreover assume that ${\rm supp }(f)$ and ${\rm supp }(g)$ are small.
Then 
$$
\|{\mathcal {T}}_{\lambda}^{r}f\,{\mathcal{T}}_{\mu}^{q}g\|_{L^2}^2
=
\int\int\int\int
K(t,s,t',s')
f(t)g(s)\overline{f(t')}\overline{g(s')}
dtdsdt'ds'
$$
where 
$$
K(t,s,t',s')=\int
e^{i\Phi_{\lambda,\mu}(x,t,t',s,s')}
a_{r}(x,t)\overline{a_{r}(x,t')}
a_{q}(x,s)\overline{a_{q}(x,s')}
dx
$$
with 
$$
\Phi_{\lambda,\mu}(x,t,t',s,s')=\lambda(\varphi_{r}(x,t)-\varphi_{r}(x,t'))+\mu(\varphi_{q}(x,s)-\varphi_{q}(x,s')).
$$
Due to Lemma \ref{le.carsjhyp},
$$
|\nabla_{x}\Phi_{\lambda,\mu}|\geq C(\lambda|t-t'|+\mu|s-s'|).
$$
Moreover
$$
|\partial^{\alpha}_{x}\Phi_{\lambda,\mu}|\leq C_{\alpha}(\lambda|t-t'|+\mu|s-s'|).
$$
Hence, by integrations by parts in $x$, we get easily, for every $N$,
$$
|K(t,s,t',s')|
\leq
C_{N}(\mu|s-s'|+\lambda|t-t'|)^{-N}.
$$
Schur's Lemma gives
$$
\|{\mathcal {T}}_{\lambda}^{r}f(x)\,{\mathcal{T}}_{\mu}^{q}g(x)\|_{L^2}^2
\leq C \left[ \int\int \frac{dtds}{(1+\mu|s|+\lambda|t|)^3}
\right]\|f\|_{L^2}^{2}\|g\|_{L^2}^{2}
$$
and therefore
$$
\|{\mathcal {T}}_{\lambda}^{r}f(x)\,{\mathcal{T}}_{\mu}^{q}g(x)\|_{L^2}
\leq 
\frac{C}{\sqrt{\mu\lambda}}\|f\|_{L^2}\|g\|_{L^2}
$$
which is better than \eqref{eq4.1}.
\\
\newline
We are left with the following two cases :
\begin{enumerate}
\item\label{casei}
$a_r$, ${a}_q$ are localized close to $t=t_0$,
$s=t_0$ respectively.
\item\label{caseii}
$a_r$, ${a}_q$ are localized close to $t=t_0$,
$s=t_0+\pi$ respectively.
\end{enumerate}
We are going to study the case \ref{casei} and will give only an outline for the case \ref{caseii} which is similar.
\\
\newline 
Let us introduce the following germs of Lagrangian manifolds in
$T^{*}(\R)$ :
\begin{equation}\begin{aligned}
\Lambda_{r}(x)&=
\big\{
\big(t,-\frac{\partial\varphi_{r}}{\partial t}(x,t)\big),
\,\,
t\sim t_0
\big\}\\
\widetilde{\Lambda}_{r}(x)&=
\big\{
\big(t,\frac{\partial\varphi_{r}}{\partial t}(x,t)\big),
\,\,
t\sim t_0 + \pi 
\big\}.
\end{aligned}\end{equation}
\begin{proposition}
\label{le4.2}
For any $q,r\in [\varepsilon/C,C\varepsilon]$, there exist two germs of canonical transformations 
$, \widetilde {\chi}_{q,r}$
\begin{gather}
 \chi_{q,r}: ( T^* \mathbb{R}, (t_0, 0)) \rightarrow  ( T^* \mathbb{R}, (t_0, 0))\\
 \widetilde{ \chi}_{q,r}: ( T^* \mathbb{R}, (t_0, 0)) \rightarrow  ( T^* \mathbb{R}, (t_0+ \pi, 0))
\end{gather}
 such that for any $x $ close to $0$,
\begin{equation}\begin{aligned}
\chi_{q,r}(\Lambda_{r}(x))&=\Lambda_{q}(x), 
\\
\widetilde{\chi}_{q,r}(\Lambda_{r}(x))&=\widetilde{\Lambda}_{q}(x), 
\end{aligned}\end{equation}
and $\chi_{q,r}$ is close to the identity whereas $\widetilde{\chi}_{q,r}$ is close to the map:
$$
(t, \tau) \mapsto (t+ \pi, \tau).
$$
\end{proposition}
This lemma will be proved in Section~\ref{subsec4.3}. Let us show how to finish the proof of Theorem~\ref{th2}. First assume that we are in case~\ref{casei}. Denote by
$$
C_{r}=
\big\{
\big(x,\frac{\partial\varphi_{r}}{\partial x},t,-\frac{\partial\varphi_{r}}{\partial t}  \big),
\,\,
t\sim t_0,\,x\sim 0
\big\}.
$$
the canonical relation associated to the Fourier integral operator ${\mathcal {T}}_{\lambda}^{r}$.
Consider  $U_{\mu}^{q,r}$ a Fourier integral operator associated with the canonical transformation $\chi_{q,r}$,
which is (locally) unitary. Then the canonical relation associated to
${\mathcal{T}}_{\mu}^{q}\circ U_{\mu}^{q,r}$ is
$$
C_{q,r}'=
\big\{
\left(x,\frac{\partial\varphi_{q}}{\partial x},\chi_{q,r}\big(s,-\frac{\partial\varphi_{q}}{\partial s}\big)\right),
\,\,
s\sim t_0,\,x\sim 0
\big\}.
$$
But we know that this canonical relation  has the form :
$$
\big\{
\left(x,\frac{\partial\psi}{\partial x},t,-\frac{\partial\psi}{\partial t}\right),
\,\,
t\sim t_0,\,x\sim 0
\big\}
$$
because $\chi_{q,r}$ is close to the identity. Taking into account the first part of Proposition~\ref{le4.2}, we get
$$
\frac{\partial \psi}{\partial t}(x,t)
=
\frac{\partial\varphi_{r}}{\partial t}(x,t)\, .
$$
Hence $\psi(x,t)=\varphi_{r}(x,t)+\theta(x)$ which implies
$$
{\mathcal{T}}_{\mu}^{q}\,\circ\,U_{\mu}^{q,r}(h)(x)
=
e^{i\mu\theta(x)}\int_{\R}
e^{i\mu\varphi_{r}(x,t)}b(x,t)h(t)dt
$$
modulo terms of order ${\cal O}(\mu^{-1})$, 
which are negligible, and where $b\in C_{0}^{\infty}(\R^2\times\R)$ satisfies the same kind of hypotheses as ${a_q}$. 
\\
\newline
We now state a bilinear version of Carleson-Sj\"olin's lemma which will be proved in section~\ref{subsec4.2}. 
\begin{proposition}\label{prop4.3}
Consider 
$\varphi \in C^\infty(\R^{2}\times\R; \R)$ and 
$a, \tilde{a} \in C^\infty_{0}(\R^{2}\times\R; \R)$ 
such that
$$
(x,t)\in {\rm supp}(a)\cup {\rm supp}(\tilde{a})
\Longrightarrow
\det(\nabla_{x}\varphi'_{t}(x,t),
\nabla_{x}\varphi''_{tt}(x,t))
\neq 0.
$$
Then if
$$
{\rm dist}\,
({\rm supp}(f),\, {\rm supp}(g))\leq \delta
$$
the operators ${\mathcal {T}}_{\lambda}$, $\widetilde{ {\mathcal {T}} }_{\mu}$
defined by 
\begin{equation}\label{eq4.7}\begin{aligned}
{\mathcal {T}}_{\lambda}f(x)&:=\int_{\R}
e^{i\lambda\varphi(x,t)}a(x,t)f(t)dt,\\
{\widetilde {\mathcal{T}}}_{\mu}f(x)&:=\int_{\R}
e^{i\mu\varphi(x,t)}\widetilde{a}(x,t)f(t)dt
\end{aligned}\end{equation}
satisfy
$$
\|{\mathcal {T}}_{\lambda}f\, \widetilde{\mathcal{T}}_{\mu}g\|_{L^2(\R^2)}
\leq
C\lambda^{-\frac{1}{4}}\,\mu^{-\frac{1}{2}}\,\|f\|_{L^2(\R)}\|g\|_{L^2(\R)}
$$
for any $1\leq\lambda\leq\mu$.
\end{proposition}
Using Proposition \ref{prop4.3}, we get
$$
\|{\mathcal {T}}_{\lambda}^{r}(f)\times ({\mathcal{T}}_{\mu}^{q}\,\circ\,U_{\mu}^{q,r})(h)\|_{L^2}
\leq
\frac{C}{\lambda^{\frac{1}{4}}\mu^{\frac{1}{2}}}
\|f\|_{L^2}\|h\|_{L^2}
$$
which is the estimate we wanted to prove (with  $g=U_{\mu}^{q,r}(h)$, since $U_{\mu}^{q,r}$ is unitary).
\\
\newline
We now give a brief outline for the case~\ref{caseii}. We use the second part of Proposition~\ref{le4.2}.  
Remark that for any $t$ close to $t_0$, we have $\widetilde {\chi} _{q,r}(t,0)= (t+\pi, 0)$.
If $\widetilde{U}_ \mu^{q,r}$ is a Fourier integral operator associated to $\widetilde{\chi}_{q,r}$, we deduce
$${{\mathcal{T}}}_\mu^q \circ \widetilde{U}_\mu^{q,r}(h) (x) = e^{i\mu \widetilde {\theta}(x)} \int _{\mathbb{R}} e^{-i\mu \varphi_r( x,t)} \widetilde{b} (x,t) h(t) dt,$$
where $\widetilde{b} \in C^\infty_0( \mathbb{R}^3)$ is localized close to $x=0,t=t_0$. Finally, to conclude, 
it is enough to remark that
$$
\|{{\mathcal{T}}}_\lambda ^r (f) \overline{ {\widetilde {\mathcal{T}}}_\mu^q(g)}\|_{L^2} 
= \|{{\mathcal{T}}}_\lambda ^r (f) { {\widetilde {\mathcal{T}}}_\mu^q(g)}\|_{L^2} 
$$
and apply Proposition~\ref{prop4.3}. 
%%%%%%%%%%%%%%%%%%%%%%%%%%%%%%%%%%%%%%%%%%%%%%%%%%%%%%%%%%%%%%%%%%%%%%%%%%%%%%%%%%%%%%%%%%%%%%%%%%%%%%%%%%%%%%%%%%%%%%%%%%%
\subsection{A bilinear Carleson-Sj\"olin Lemma}\label{subsec4.2}
We now prove Proposition~\ref{prop4.3}. 
According to \eqref{eq4.7} and the compactness of the supports of $a$, $\widetilde{a}$,
we can reduce the study to the case where $a$, $\widetilde{a}$ are  
supported in $K\times I$, with $K\subset \subset {\mathbb R}^2$ and  $I$ is a small interval in $\R$ which can, by translation invariance, be supposed to be equal to $[-\frac{\delta}{2},\frac{\delta}{2}]$. Then we have
$$
{\mathcal {T}}_{\lambda}f(x)\, \widetilde{\mathcal {T}}_{\mu}g(x)
=
\int_{I}\int_{I}
e^{i\mu\Phi_{\epsilon}(x,t,s)}a(x,t)\,\widetilde{a}(x,s)f(t)g(s)dt ds
$$
where $\Phi_{\epsilon}(x,t,s)=\epsilon\varphi(x,t)+\varphi(x,s)$, $\epsilon=\frac{\lambda}{\mu}\leq 1$.\par
For $h\in L^{2}(\R\times\R)$ denote by
$$
S^{\pm}_{\lambda,\mu}h(x)=\int_{(I\times I)\cap \{\pm(t-s)>0\}}e^{i\mu\Phi_{\epsilon}(x,t,s)}a(x,t)\,\widetilde {a}(x,s)h(t,s)dt ds
$$
We have ${\mathcal {T}}_{\lambda}f\, {\mathcal {T}}_{\mu}g= S^{+}_{\lambda,\mu}(f\otimes g)+S^{-}_{\lambda,\mu}(f\otimes g)$ and to prove the lemma, it is enough to show that
$$
\|S^{\pm}_{\lambda,\mu}\|_{L^{2}(\R)\rightarrow L^{2}(\R)}
\leq C\lambda^{-\frac{1}{4}}\,\mu^{-\frac{1}{2}}.
$$
Let us study $S^{+}_{\lambda,\mu}$ (the case of $S^{-}_{\lambda,\mu}$ is similar). Compute
\begin{equation*}
\|S^{+}_{\lambda,\mu}h\|_{L^2(\R^2)}^{2}=
\int_{((I\times I)\cap \{t>s\})^2}
K_{\lambda,\mu}(t,s,t',s')h(t,s)\overline{h(t',s')}dtdsdt'ds'
\end{equation*}
where
\begin{equation}\label{eq4.20}
K_{\lambda,\mu}(t,s,t',s')=
\int_{\R^2}e^{i\mu\left[\Phi_{\epsilon}(x,t,s)-\Phi_{\epsilon}(x,t',s')\right]}
a(x,t)\widetilde {a}(x,s)\overline{a(x,t')}\,\,\overline{\widetilde {a}(x,s')}\,dx.
\end{equation}
On  $(I\times I)\cap \{t>s\}$, 
let us perform the bijective change of variables
from $\{t>s\}$ to  $\{(u,v)\,: \,u>0\}$
\begin{equation*}
\left\{
\begin{array}{l}
u  =\frac{\epsilon}{2}(t-s)^2 \\
v  = s+\epsilon t
\end{array}
\right.
\end{equation*}
whose Jacobian is
$$
\frac{D(u,v)}{D(t,s)}=
\left|
\begin{array}{rl}
\epsilon(t-s) & \epsilon
\\
-\epsilon(t-s) & 1
\end{array}
\right|
=
\epsilon(t-s)(1+\epsilon)=(1+\epsilon)\sqrt{2\epsilon u}.
$$
Set $\Phi_{\epsilon}(x,t,s)=\widetilde{\Phi}_{\epsilon}(x,u,v)$ for any $(u,v)$ in the image 
$\Omega$ of
$(I\times I)\cap \{t>s\}$. Let us show that $\widetilde{\Phi}_{\epsilon}$ is a uniformly non degenerated phase on
$K\times \Omega$. Compute
\begin{equation*}\begin{aligned}
\Delta_{\epsilon}:=
\left|
\begin{array}{ll}
\widetilde{\Phi}_{\epsilon,\, x_1 u}^{''} &  \widetilde{\Phi}_{\epsilon,\, x_2 u}^{''}
\\
\widetilde{\Phi}_{\epsilon,\, x_1 v}^{''} &  \widetilde{\Phi}_{\epsilon,\, x_2 v}^{''}
\end{array}
\right|
& = 
\frac{D(t,s)}{D(u,v)}
\left|
\begin{array}{ll}
{\Phi}_{\epsilon,\, x_1 t}^{''} &  {\Phi}_{\epsilon,\, x_2 t}^{''}
\\
{\Phi}_{\epsilon,\, x_1 s}^{''} &  {\Phi}_{\epsilon,\, x_2 s}^{''}
\end{array}
\right|
\\
& = 
\frac{D(t,s)}{D(u,v)}
\left|
\begin{array}{ll}
\epsilon {\varphi}_{\epsilon,\, x_1 t}^{''}(x,t) &  \epsilon{\varphi}_{\epsilon,\, x_2 t}^{''}(x,t)
\\
{\varphi}_{\epsilon,\, x_1 s}^{''}(x,s) &  {\varphi}_{\epsilon,\, x_2 s}^{''}(x,s)
\end{array}
\right|
\\
& = 
\frac{D(t,s)}{D(u,v)}
\left[
(s-t)\epsilon
\left|
\begin{array}{ll}
{\varphi}_{\epsilon,\, x_1 t}^{''}(x,t) &  {\varphi}_{\epsilon,\, x_2 t}^{''}(x,t)
\\
{\varphi}_{\epsilon,\, x_1 tt}^{'''}(x,t) &  {\varphi}_{\epsilon,\, x_2 tt}^{'''}(x,t)
\end{array}
\right|
+
{\mathcal O}((s-t)^2)
\right]
\end{aligned}
\end{equation*}
and if $\delta>0$ is small enough we get :
$$
0<\alpha\leq |\Delta_{\epsilon}|\leq \beta.
$$
On the other hand, since
\begin{equation*}\begin{aligned}
\frac{\partial}{\partial t}& =  \epsilon(t-s) \frac{\partial}{\partial u}+\epsilon \frac{\partial}{\partial v}
\\
\frac{\partial}{\partial s} & = \epsilon(s-t) \frac{\partial}{\partial u}+\frac{\partial}{\partial v}
\end{aligned}
\end{equation*}
we have
\begin{equation*}\begin{aligned}
\frac{\partial }{\partial u} & = \frac{1}{(1+\epsilon)\epsilon(t-s)} 
\left(\frac{\partial}{\partial t}-\epsilon \frac{\partial}{\partial s}\right)
\\
\frac{\partial}{\partial v} & =  \frac{1}{(1+\epsilon)} 
\left(\frac{\partial}{\partial t}+\frac{\partial}{\partial s}\right)
\end{aligned}
\end{equation*}
which implies
\begin{equation*}\begin{aligned}
\frac{\partial\widetilde{\Phi}_{\epsilon}}{\partial u} & =  \frac{1}{(1+\epsilon)\epsilon(t-s)} 
\left[\epsilon\frac{\partial\varphi}{\partial t}(x,t)-\epsilon \frac{\partial\varphi}{\partial s}(x,s)\right]
\\
\frac{\partial\widetilde{\Phi}_{\epsilon}}
{\partial v} & =  \frac{1}{(1+\epsilon)} 
\left[\epsilon\frac{\partial \varphi}{\partial t}(x,t)+\frac{\partial\varphi}{\partial s}(x,s)\right].
\end{aligned}\end{equation*}
We deduce that
$
\frac{\partial^{2}\widetilde{\Phi}_{\epsilon}}
{\partial x\partial u},
$
$
\frac{\partial^{2}\widetilde{\Phi}_{\epsilon}}
{\partial x\partial v}
$
are bounded on $K\times \Omega$.
Since $\Omega$ is included in a convex set where these properties still hold, we obtain that there exists a constant $C>0$ (independent of $\epsilon >0$) such that
$$
C\left(|u-u'|+|v-v'|\right)
\leq
|\nabla_{x}\widetilde{\Phi}_{\epsilon}(x,u,v)-\nabla_{x}\widetilde{\Phi}_{\epsilon}(x,u',v')|
\leq\frac{1}{C}\left(|u-u'|+|v-v'|\right).
$$
and more generally 
$$|\partial^\alpha_{x} \widetilde \Phi_{\varepsilon} (x,u,v)-\partial^\alpha_{x} \widetilde \Phi_{\varepsilon} (x,u',v')|\leq C_{\alpha}\left( |u-u'|+|v-v'|\right)$$
Consequently, by integrations by parts with respect to $x$ in~\eqref{eq4.20}, we obtain for any $N\in\mathbb{N}$
$$
|K_{\lambda,\mu}(t,s,t',s')|
\leq C_N
(1+\mu|u-u'|+\mu|v-v'|)^{-N}
$$
for any $(t,s)$, $(t',s')$ in $(I\times I)\cap \{t>s\}$.\par
To conclude, we apply Schur's lemma: since $K_{\lambda,\mu}$ is Hermitian,
\begin{equation*}\begin{aligned}
\|S^{+}_{\lambda,\mu}\|^{2}_{L^2\rightarrow L^2}
& \leq  \sup_{(t,s)}\int_{\R^2}|K_{\lambda,\mu}(t,s,t',s')|dt'ds'
\\
& \leq  C\sup_{(u,v)\in\Omega}\,\,
\int_{(I\times I)\cap \{t>s\}}(1+\mu|u-u'|+\mu|v-v'|)^{-4}dt'ds'
\\
& \leq 
C \sup_{(u,v)\in\R^2}
\int_{\R^2}(1+\mu|u-u'|+\mu|v-v'|)^{-4}\frac{du'dv'}{\sqrt{2\epsilon|u'|}}.
\end{aligned}\end{equation*}
 Set $u'=u+\frac{z}{\mu}$,
$v'=v+\frac{w}{\mu}$
\begin{equation*}
\begin{aligned}
\|S^{+}_{\lambda,\mu}\|^{2}_{L^2\rightarrow L^2}
& \leq 
\frac{C}{\mu^2}\sqrt{\frac{\mu}{2\epsilon}}\,\,
\sup_{\tilde{u}\in\R}
\int_{\R^2}(1+|z|+|w|)^{-4}\frac{dzdw}{\sqrt{|z+\tilde{u}|}}
\\
& \leq 
\frac{\widetilde{C}}{\lambda^{\frac{1}{2}}\mu}
\end{aligned}
\end{equation*}
since $\epsilon=\frac{\lambda}{\mu}$.
\qed 
%%%%%%%%%%%%%%%%%%%%%%%%%%%%%%%%%%%%%%%%%%%%%%%%%%%%%%%%%%%%%%%%%%%%%%%%%%%%%%%%%%%%%%%%%%%%%%%%%%%%%%%%%%%%%%%%%%%%%%%%%%%%%%%%
\subsection{Canonical Transformations}\label{subsec4.3}
We are going to prove Proposition~\ref{le4.2}. Let us first deal with $\chi_{q,r}$. We have to define a germ of canonical transformation
$$
(t,\tau)\mapsto(s,\sigma)=(S_{q,r}(t,\tau),\Sigma_{q,r}(t,\tau))
$$
such that for $t\sim t_0$, $x\sim 0$,
\begin{equation}\label{eq4.8}
-\frac{\partial\varphi_{q}}{\partial s}
\left(x,S_{q,r}\big(t,-\frac{\partial\varphi_{r}}{\partial t}(x,t)\big)
\right)=\Sigma_{q,r}\big(t,-\frac{\partial\varphi_{r}}{\partial t}(x,t)\big).
\end{equation}
It is natural to study, for $\tau\in\R$ close to $0$
and $t$ close to  $t_0$, the structure of the set
$$
\left\{
x\in M,x\sim 0,-\frac{\partial\varphi_{r}}{\partial t}(x,t)=\tau
\right\}.
$$
Since $\nabla_{x}\partial_{t}\varphi_{r}\neq 0$, we know that this set is a smooth curve close to $0$. 
We are going to show that this curve is a geodesic which 
passes through $\exp_{o}(re^{it})$.
Let $x$ be a point on this curve.
Let  $u=u_{r}(x,t)\in T_{x}M$ be the initial speed of the geodesic which joins $x$ to $\exp_{o}(re^{it})$.
Define for any $r,t$,
$$
H(r,t)=
\exp_{o}(re^{it}).
$$
We have
$$
H(r,t)=
\exp_{x}(-\varphi_{r}(x,t)u_{r}(x,t)).
$$
Differentiating this identity with respect to $t$ we get
$$
\partial_{t}H(r,t)
=
T_{-\varphi_{r}u_{r}}(\exp_{x})
\left[
-\partial_{t}\varphi_{r}u_{r}-\varphi_{r}\partial_{t}u_{r}
\right].
$$
Denote by $v_{r}(x,t)=T_{-\varphi_{r}u_{r}}(\exp_{x})(u_r)\in T_{H(r,t)}(M)$. Taking the scalar product with  $v_r$,
we get, using Gauss' lemma,
\begin{equation}\label{eq4.9}
-\partial_{t}\varphi_{r}=
g(v_r,\partial_t H_r),
\end{equation}
so that $-\partial_t \varphi_r$ depends upon $x$ only through $v_r$, 
which is the speed at $H(t,r)$ of the geodesic joining $x$ to $H(r,t)$. 
As a consequence, $-\partial_t \varphi_r$ is constant on this geodesic.
\begin{figure}
$$\ecriture{\includegraphics[width=8cm]{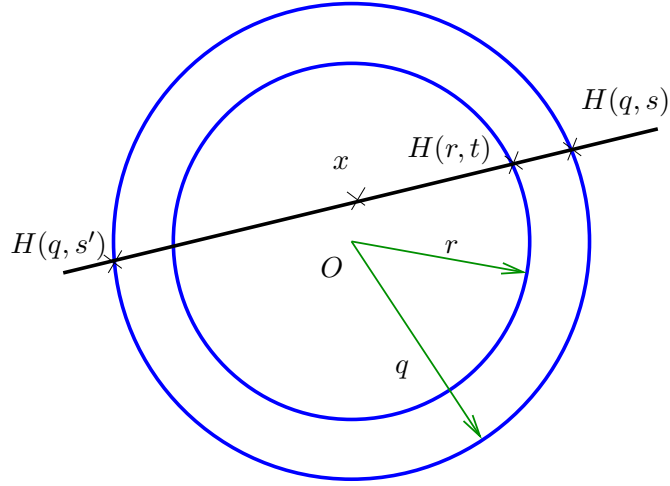}}
{\aat{-3}{19}{$H(q, s')$}\aat{23}{26}{$x$}\aat{29}{27}{$H(r,t)$}\aat{43}{31}{$H(q,s)$}\aat{32}{19}{$r$}\aat{28}{9}{$q$}\aat{22}{17}{$O$}}$$
\caption{The canonical transformations $\chi_{q,r}$ et $\widetilde{\chi} _{q,r}$ }
\end{figure}
\\
\newline
Since geodesic circles are transversal to geodesics, this geodesic intersects the circle of radius $q$ at 
two points $H(q,s)$ and $H(q,s')$. 
Since $q$ is close to $r$, we can distinguish between these points by asking that  $s$ be close to $t_0$, and that $s'$ be close 
to  $t_0+\pi$.
Let  $S_{q,r}$ be defined by 
$$
s=S_{q,r}(t,\tau)
$$
Furthermore, since  $H(q,s)$ is on the geodesic, the quantity $-\partial_{s}\varphi_{q}(x,s)$ 
depends only on $(t,\tau)$. We denote it by 
$$
\sigma=\Sigma_{q,r}(t,\tau).
$$
All that remains to do is proving that the map $(S_{q,r},\Sigma_{q,r})$ is canonical.
The vector field $v_r=v_r(t,\tau)$ above $H(r,t)$ is completely defined by the combination of~\eqref{eq4.9}, the fact that it is a unitary vector 
and the knowledge of its value $\partial_{r}H(r,t)$ at $\tau=0$.
\begin{lemme}\label{le4.3} 
For any $r\geq \epsilon/C$, the map
$$\kappa_r: (t, \tau) \in T^* \mathbb{R} \mapsto (H(r,t), v_r(t, \tau))\in TM$$
defined for $(t,\tau)$ close to $ (t_0,0)$ is a symplectic embedding with values in
$$S_r(M)= \{ (y,v)\in TM; d(y,0) = r, g_y(v,v)=1, g_y(\nabla d(y,0),v)\neq 0\},$$
where $TM$ is equipped with the symplectic structure of $T^*M$ inherited from $g$.
\end{lemme}
\proof We work in geodesic coordinates $y= H(r,t)$. The vector $v_r(t, \tau)$ corresponds, via $g$, to the covector
$$ \xi_r= \rho dr + \theta dt$$
so that the equation~\eqref{eq4.9} is
$$\tau = g(v_r, \partial_t H) = \theta$$
Since $r$ is constant on $S_r(M)$, the restriction to $S_r(M)$ of the symplectic form 
$d\rho \wedge dr + d\theta \wedge dt$ is exactly $d\tau \wedge dt$. 
Finally the fact that $S_r(M)$ is symplectic is ensured by the condition $g_y( \nabla_y d(y,0),v) \neq 0$ 
which implies the non vanishing of the Poisson bracket of $g(v,v) $ and $d(y,0)$.
\endproof
The following lemma is standard:
\begin{lemme}\label{le4.4} 
Let $S$ be a symplectic manifold and $p$ a function with real values on $S$ such that $dp\neq 0$. 
Let $S_1, S_2$ be two symplectic submanifolds of $S$ of codimension $2$, included in $\{p=0\}$. 
Let $\rho_1\in S_1$, $\rho_2\in S_2$ and $T\in \mathbb{R}$ be such that $\exp(TH_p) (\rho_1) = \rho_2$. 
Then for any $\rho$ close to $\rho_1$ there exists a unique $T(\rho)$ close to $T$, such that 
$\exp(T(\rho)H_p)(\rho)= F(\rho) \in S_2$. Furthermore the map
$$F:S_1\longrightarrow S_2$$ is canonical.
\end{lemme}
\proof 
Composing by $\text{exp}(TH_p)$, we reduce the problem to the case where $T=0$, $\rho_1= \rho_2$ and the problem is local. We choose symplectic coordinates $(x, \xi)$ such that $p= \xi_1$. $S_1$ and $S_2$ are given by $\{\xi_1= f_1=0\}$ and $\{\xi_1= f_2=0\}$ respectively, with $\partial f_j/\partial x_1\neq 0$. Consequently we can replace $f_j$ by $x_j - g_j(x', \xi')$. Then
$$\text{exp}( t H_p)( g_1(x',\xi'), x', 0, \xi')= (g_1(x',\xi')+ t, x', 0, \xi')$$
 $$T(x',\xi')= g_2(x',\xi')- g_1(x',\xi')$$
and the map $F$ is the identity in the coordinate system $(x', \xi ')$.
\endproof
We take $p(m,v)= g_m(v,v)-1$ and in case~\ref{casei}, 
\begin{equation}
\begin{gathered}
S_1= S_r(M), \qquad S_2= S_q(M)\\
\rho_1= (H(r, t_0), \partial_rH(r, t_0)) \text{ and } \rho_2= (H(q, t_0), \partial_rH(q, t_0)).
\end{gathered}
\end{equation}
 This choice gives a canonical transformation $F_{q,r}$ and we check easily that with the notations of Lemma~\ref{le4.3},
$$ \chi_{q,r}= \kappa_q^{-1}\circ  F_{q,r} \circ \kappa _r.$$
For case~\ref{caseii}, we apply Lemma~\ref{le4.4} with
$S_1= S_r(M)$, $S_2= S_q(M)$, $\rho_1= (H(r, t_0), \partial_rH(r, t_0))$ and 
$\rho_2= (H(q, t_0+\pi), -\partial_qH(q, t_0+\pi))$. 
This choice gives a canonical transformation $\widetilde {F}_{q,r}$ and we check that
$$ 
\widetilde{\chi}_{q,r}= \kappa_q^{-1}\circ  \widetilde{F}_{q,r} \circ \kappa _r
$$
satisfies, close to $0$,
$$ \widetilde{\chi}_{q,r}( \Lambda_r(x))= \{ (s, \frac{\partial \varphi_q} { \partial s}), s \sim t_0 + \pi\}.$$
\endproof 
%\bibliographystyle{plain}
%\bibliography{/Users/nicolasb/biblio/biblio}
%\bibliography{/home/u1/vis/burq/biblio/biblio}

\def\cprime{$'$}

\end{document}